\documentclass[11pt,a4paper]{article}
\usepackage{amsmath,amssymb,subcaption,url}
\usepackage{fullpage}
\usepackage{tikz}
\usetikzlibrary{shapes}

\newcommand{\mmw}{\mathit{MM}}
\newcommand{\mme}{\mathit{MM}_{\textnormal{e}}}
\newcommand{\abw}{\mathit{AB}}
\newcommand{\abe}{\mathit{AB}_{\textnormal{e}}}

\newcommand{\ans}[2]{{#1}\mbox{B}{#2}\mbox{W}}

\newcommand{\ac}{\mbox{a}}

\newcommand{\opt}[1]{\textbf{#1}}

\begin{document}

\title{Algorithm to Prove Formulas for the Expected Number of Questions in Mastermind Games}
\author{Marcin Peczarski\\
\\
Institute of Informatics, University of Warsaw\\
ul.\ Banacha 2, PL-02-097 Warszawa, Poland\\
email: marpe@mimuw.edu.pl}
\maketitle

\begin{abstract}
\noindent
We close the gap in the proof (published by Chen and Lin) of formulas for the
minimum number of questions required in the expected case for Mastermind and its
variant called AB game, where both games are played with two pegs and $n$
colors.
For this purpose, we introduce a new model to represent the game guessing
process and we develop an algorithm with automatizes the proof.
In contrary to the model used by Chen and Lin, called graph-partition approach,
which is limited to two pegs, our model and algorithm are parametrized with the
number of pegs and they could potentially be used for any number of pegs.
\end{abstract}



\section{Introduction}

Generalized Mastermind is a deductive game for two players:
a \emph{codemaker} and a \emph{codebreaker}.
The codemaker chooses a secret consisting of $p$ pegs out of $n$ possible
colors.
We denote colors using non-negative numbers: $0,1,2,\dots,n-1$.
We denote a secret by $(s_1,s_2,s_3,\dots,s_p)$, where peg $i$ has the color
$s_i\in\{0,1,2,\dots,n-1\}$ for $i=1,2,3,\dots,p$.
The codebreaker tries to guess the secret by asking questions.
We denote a question by $(q_1,q_2,q_3,\dots,q_p)$, where
$q_i\in\{0,1,2,\dots,n-1\}$ is the guessed peg color in position $i$ for
$i=1,2,3,\dots,p$.
The codemaker responds to a question with $b$ black and $w$ white pegs.
A black peg means that a peg in the question is correct in both position and
color.
A white peg means that a peg in the question is correct in color but not in
position.
The formal definition is as follows:
{\setlength\arraycolsep{2pt}
\begin{eqnarray*}
b&=&|\{i:s_i=q_i\}|,\\
w&=&-b+\sum_{j=0}^{n-1}\min\{|\{i:s_i=j\}|,|\{i:q_i=j\}|\}.
\end{eqnarray*}}%
Further, we denote an answer by $\ans{b}{w}$.
The game ends when the codemaker answers with $\ans{p}{0}$.
Original Mastermind is a famous game invented by Mordecai Meirowitz in 1970,
and is played with $p=4$ pegs and $n=6$ colors.
We omit the adjective ``generalized'', having in mind the game with $p$ pegs and
$n$ colors.

In Mastermind, repeated colors are allowed in secrets and questions.
Another deductive game is a variant of Mastermind called ``Bulls an Cows'' or AB
game.
The difference is that AB game allows only secrets and questions without
repeated colors.
Hence, it must hold that $n \ge p$.
The original AB game is played with $p=4$ pegs and $n=10$ colors (digits).

There are two ways of optimizing the codebreaker's strategy, namely minimizing
the number of questions in the worst case or in the expected case.
We denote the minimal number of the questions required by the codebreaker in the
worst or expected case by $\mmw(p,n)$, $\mme(p,n)$ for Mastermind and
$\abw(p,n)$, $\abe(p,n)$ for AB game, respectively.
The last question, which is answered with $p$ black pegs, is also counted.

The value of $\mmw(4,6)=5$ is found by Knuth \cite{K76}.
The value of $\mme(4,6)=5625/1296\approx4.34$ is computed by Koyama and Lai
\cite{KL93}.
Some other values of $\mmw(p,n)$ and $\mme(p,n)$ for small $p$ and $n$ are
found by Goddard \cite{G04}.
Note that in the table of \cite{G04} the values for $\mmw(2,n)$ are shifted by
one color.
A few more values of $\mmw(p,n)$ are computed by J\"ager and Peczarski
\cite{JP09}.
The value of $\mme(4,7)=11228/2401\approx4.676$ is reported by Ville \cite{V13}.
The value of $\abw(4,10)=7$ is computed by Chen et al.\ \cite{ChLHH07}.
The value of $\abe(4,10)=26274/5040\approx5.213$ is reported by Huang et
al.\ \cite{HChL09}.

The formula $\mmw(2,n)=\lfloor n/2 \rfloor + 2$ for $n \ge 2$ is
independently given in \cite{ChL04mm} and proved in \cite{G04}, and later also
proved in \cite{JP09}.
The formula $\mmw(3,n)=\lfloor (n-1)/3 \rfloor + 4$ for $n \ge 5$ and
tight bound $4 \le \mmw(4,n) - \lfloor n/4 \rfloor \le 6$ for $n \ge 16$
are proved in \cite{JP09}.
The asymptotic formula $\mme(2,n) = n/3 + 17/8 + o(1)$ is proved in \cite{G04}.
The following formula (true for $n \ge 3$) is given in \cite{ChL04mm}:
\begin{equation}\label{eq:mme}
\mme(2,n) = \left\{
  \begin{array}{l l}
    (8n^3 + 51n^2 - 74n + 48)/24n^2, & \textnormal{if $n$ is even,}\\
    (8n^3 + 51n^2 - 80n + 69)/24n^2, & \textnormal{if $n$ is odd.}
  \end{array} \right.
\end{equation}
The following formula is given in \cite{ChL04ab} and proved in \cite{JP15}:
\begin{displaymath}
  \abw(2,n) = \lceil n/2 \rceil + 1
    = \lfloor (n+1)/2 \rfloor + 1, \quad\textnormal{if $n \ge 2$}.
\end{displaymath}
It is also proved in \cite{JP15} that
\begin{displaymath}
  \begin{array}{ll}
    \abw(3,n) = \lfloor n/3 \rfloor + 3,     & \textnormal{if $3 \le n \le 7$,}\\
    \abw(3,n) = \lfloor (n+1)/3 \rfloor + 3, & \textnormal{if $n \ge 8$,}\\
    \abw(4,n) = \lfloor (n+2)/3 \rfloor + 3, & \textnormal{if $4 \le n \le 11$},\\
    \abw(4,n) = 8                          , & \textnormal{if $n = 12,13$},\\
    \abw(4,n) \ge \lfloor (n+3)/4 \rfloor + 4, & \textnormal{if $n \ge 14$},\\
    \abw(4,n) \le \lfloor (n+3)/4 \rfloor + 5, & \textnormal{if $n \ge 14$}.
  \end{array}
\end{displaymath}
The following formula (true for $n \ge 2$) is given in \cite{ChL04ab}:
\begin{equation}\label{eq:abe}
\abe(2,n) = \left\{
  \begin{array}{l l}
    (4n^3 + 21n^2 - 76n + 72)/12n(n-1),  & \textnormal{if $n$ is even,}\\
    (4n^3 + 21n^2 - 82n + 105)/12n(n-1), & \textnormal{if $n$ is odd.}
  \end{array} \right.
\end{equation}

We write above that the formulas for $p=2$ pegs in \cite{ChL04mm,ChL04ab} are
``given'' and not ``proved'', because proofs in these papers have a gap.
The authors assume that only colors already possible for the secret are allowed
for a question.
For example, if $(0,0)$ is the first question asked by the codebreaker in the
game with $n$ colors and the codemaker answers with $\ans{0}{0}$ then color $0$
is not allowed in the secret in any position and other $n-1$ colors are allowed
in both positions.
The authors claim that the game is then equivalent to the game with $n-1$
colors, which is false, because (according to the game rules) the codebreaker
should be still allowed to ask questions containing color $0$.
As shown in \cite{JP11} for the Mastermind variant without white pegs in answers
played with $p=3,4,5,6,7,8$ pegs and $n=2$ colors, an additional color (allowed
in questions, but not in the secret) decrease the number of required questions
in the worst case.

The first reason of this paper is to close the gap in proofs.
Because the worst case for two pegs Mastermind and AB game, i.e.\ the formulas
for $\mmw(2,n)$ and $\abw(2,n)$, are already proved in \cite{G04,JP09,JP15},
we focus on proving the formulas for the expected case, i.e.\ equations
(\ref{eq:mme}) and (\ref{eq:abe}).
The another reason is to present a new model for the guessing process in the
Mastermind games and an algorithm which automatizes the proof.

The paper is organized as follows.
We introduce the game model in Section \ref{sec:model}.
We describe the algorithm in Section \ref{sec:algo}.
Section \ref{sec:mm} contains results of computations for Mastermind and
analysis of these results.
Section \ref{sec:ab} contains analogous results for AB game.

\section{Game Model}\label{sec:model}

The Mastermind game guessing process is represented as a tree.
A node symbolizes a game state which is a set of possible secrets, further
called a \emph{maset}.
In other words, a maset contains actual knowledge of the codebreaker about the
secret.
An internal node (not a leaf) contains also a question asked in that game state.
An edge symbolizes an answer.
A leaf is a node to which leads the edge with the answer $\ans{p}{0}$.

We represent a maset as a regular expression, further called a
\emph{maset pattern}.
A maset pattern is an alternative of \emph{clauses}.
Each clause is a \emph{tuple} of $p$ elements.
These elements could be \emph{explicit colors}, i.e.\ the numbers 0, 1, 2,
\dots, $n-1$, or a \emph{star}.
The star represents \emph{implicit colors}, which are not explicitly listed in
clauses.
The star has an index which is the number of represented colors.
All stars in a maset pattern have the same index.
A clause represents a set of secrets.
Sets of secrets represented by clauses in a maset pattern are disjoint.
For example, the maset pattern $(*_n,*_n)$ represents all secrets in the game
with $p=2$ pegs and $n$ colors.
Another example is the maset pattern
$(0,0) \mid (0,*_{n-2}) \mid (1,1) \mid (*_{n-2},1)$.
The clause $(0,*_{n-2})$ represents secrets $(0,2)$, $(0,3)$, $(0,4)$, \dots,
$(0,n-1)$.
The clause $(*_{n-2},1)$ represents secrets $(2,1)$, $(3,1)$, $(4,1)$, \dots,
$(n-1,1)$.
The empty maset pattern does not contain any clause, it does not represent any
secret and it is denoted by $\emptyset$.

A maset pattern is \emph{normalized}, if the set of all explicit colors in its
clauses contains the consecutive numbers starting from 0.
For example, the maset pattern $(1,*_{n-2})$ is not normalized and it represents
secrets $(1,2)$, $(1,3)$, $(1,4)$, \dots, $(1,n-1)$.
After normalization (color remuneration), it is equivalent to the maset pattern
$(0,*_{n-2})$ which represents secrets $(0,1)$, $(0,2)$, $(0,3)$, \dots,
$(0,n-2)$.

We need to recognize \emph{isomorphic} maset patterns.
Let $m$ be the number of clauses in the maset pattern.
We represent the maset pattern as a table with $m$ rows and $p$ columns.
An element in $i$th row and $j$th column is $j$th element of $i$th clause.
We ignore star indices.
Two maset patterns are isomorphic if there exists a permutation of rows,
a permutation of columns and a permutation of explicit colors mapping one table
to another.
For example, maset patterns $(0,*_{n-2})$ and $(*_{n-2},0)$ are isomorphic.
Other examples of isomorphic maset patterns are:
$(*_n,*_n)$ and $(*_{n-1},*_{n-1})$,
$(0,*_n)$ and $(*_{n-1},0)$,
$(0,0) \mid (*_{n-2},0)$ and $(0,0) \mid (0,*_{n-2})$,
$(0,*_{n-2}) \mid (*_{n-2},1)$ and $(1,*_{n-2}) \mid (*_{n-2},0)$.

In the algorithm presented in the next section, we use an operation which
\emph{extend} a normalized maset pattern with $u$ explicit colors to a maset
pattern with $v$ explicit colors, where $v>u$.
Extending maset pattern replaces each star $*_{n-t}$ by colors $u$,
$u+1$, $u+2$, \dots, $v-1$ and the star $*_{n-r}$, where $r=t-u+v$.
For example, if we extend the maset pattern $(*_n,*_n)$ with $u=0$ explicit
colors to $v=1$ explicit colors, we get the maset pattern
$(0,0) \mid (0,1) \mid (1,0) \mid (1,1) \mid (*_{n-1},*_{n-1})$.
If we extend the maset pattern
$(0,0) \mid (0,*_{n-2}) \mid (1,1) \mid (*_{n-2},1)$ with $u=2$ explicit colors
to $v=4$ explicit colors, we get the maset pattern
$(0,0) \mid (0,2) \mid (0,3) \mid (0,*_{n-4}) \mid
 (1,1) \mid (2,1) \mid (3,1) \mid (*_{n-4},1)$.

To close the gap in proofs of equations (\ref{eq:mme}) and (\ref{eq:abe}),
we consider a game with an \emph{additional color}, i.e.\ the codebreaker is
allowed to ask questions containing a color which is not a legal color for the
secret.
We will show that using the additional color does not help, i.e.\ does not
decrease the number of questions.
Hence, we will prove that it is sufficient for the optimal codebreaker's
strategy to use questions which contain only colors that are already not
excluded to being in the secret, as assumed in \cite{ChL04mm,ChL04ab}.
We denote the additional color by letter ``$\ac$''.
We never use the question which contains only the additional color, i.e.\ the
question $(\ac,\ac,\ac,\dots,\ac)$, as it does not distinguish secrets.

We need to generate \emph{non-isomorphic questions} respect to a given maset
pattern.
As previously, $m$ is the number of clauses in the maset pattern.
For a question we built a table with $m+1$ rows and $p$ columns.
The first $m$ rows represent the maset pattern, as previously.
The $(m+1)$th row represents the question.
Two questions are isomorphic respect to the maset pattern if there exists a
permutation of rows, a permutation of columns and a permutation of explicit
colors mapping one table to another.
For example, there are three questions pairwise non-isomorphic respect to the
maset pattern $(*_n,*_n)$, namely questions $(0,0)$, $(0,1)$, $(0,\ac)$;
there are ten questions pairwise non-isomorphic respect to the maset pattern
$(0,0) \mid (0,*_{n-2}) \mid (1,1) \mid (*_{n-2},1)$, namely questions
$(0,0)$, $(0,1)$, $(0,2)$, $(0,\ac)$, $(1,0)$, $(1,2)$, $(1,\ac)$, $(2,2)$,
$(2,3)$, $(2,a)$.

\begin{figure}[t]
  \begin{center}
    \begin{tikzpicture}
      \tikzstyle{answer}=[fill=white, inner sep=0.04cm]
      \tikzstyle{internal}=[draw, rectangle split, rectangle split parts=2]
      \tikzstyle{leaf}=[draw, ellipse, inner sep=0.04cm]
      \tikzstyle{subtree}=[]
      \node[internal] (m0) at (5,10) {$(*_n,*_n)$ \nodepart{second} $(0,1)$};
      \node[internal] (m11) at (0,6.5) {$(*_{n-2},*_{n-2})$ \nodepart{second} $(0,1)$};
      \node[internal] (m12) at (3,5) {$(1,*_{n-2}) \mid (*_{n-2},0)$ \nodepart{second} $(1,2)$};
      \node[internal] (m13) at (7.8,3.5) {$(0,0) \mid (0,*_{n-2}) \mid (1,1) \mid (*_{n-2},1)$ \nodepart{second} $(1,2)$};
      \node[internal] (m14) at (10,7.5) {$(1,0)$ \nodepart{second} $(1,0)$};
      \node[leaf] (m15) at (10,9) {$(0,1)$};
      \node[leaf] (m24) at (10,5.2) {$(1,0)$};
      \node[subtree] (s111) at (-1,5) {};
      \node[subtree] (s112) at (-0.5,5) {};
      \node[subtree] (s113) at (0,5) {};
      \node[subtree] (s114) at (0.5,5) {};
      \node[subtree] (s115) at (1,5) {};
      \node[subtree] (s121) at (2.25,3.5) {};
      \node[subtree] (s122) at (2.75,3.5) {};
      \node[subtree] (s123) at (3.25,3.5) {};
      \node[subtree] (s124) at (3.75,3.5) {};
      \node[subtree] (s131) at (7.05,2) {};
      \node[subtree] (s132) at (7.55,2) {};
      \node[subtree] (s133) at (8.05,2) {};
      \node[subtree] (s134) at (8.55,2) {};
      \draw[-latex] (m0) -- node[answer] {$\ans{0}{0}$} (m11);
      \draw[-latex] (m0) -- node[answer] {$\ans{0}{1}$} (m12);
      \draw[-latex] (m0) -- node[answer] {$\ans{1}{0}$} (m13);
      \draw[-latex] (m0) -- node[answer] {$\ans{0}{2}$} (m14);
      \draw[-latex] (m0) -- node[answer] {$\ans{2}{0}$} (m15);
      \draw[-latex] (m14) -- node[answer] {$\ans{2}{0}$} (m24);
      \draw[-latex] (m11) -- (s111);
      \draw[-latex] (m11) -- (s112);
      \draw[-latex] (m11) -- (s113);
      \draw[-latex] (m11) -- (s114);
      \draw[-latex] (m11) -- (s115);
      \draw[-latex] (m12) -- (s121);
      \draw[-latex] (m12) -- (s122);
      \draw[-latex] (m12) -- (s123);
      \draw[-latex] (m12) -- (s124);
      \draw[-latex] (m13) -- (s131);
      \draw[-latex] (m13) -- (s132);
      \draw[-latex] (m13) -- (s133);
      \draw[-latex] (m13) -- (s134);
    \end{tikzpicture}
  \end{center}
  \caption{A fragment of the Mastermind game tree for two pegs}\label{fig:three}
\end{figure}
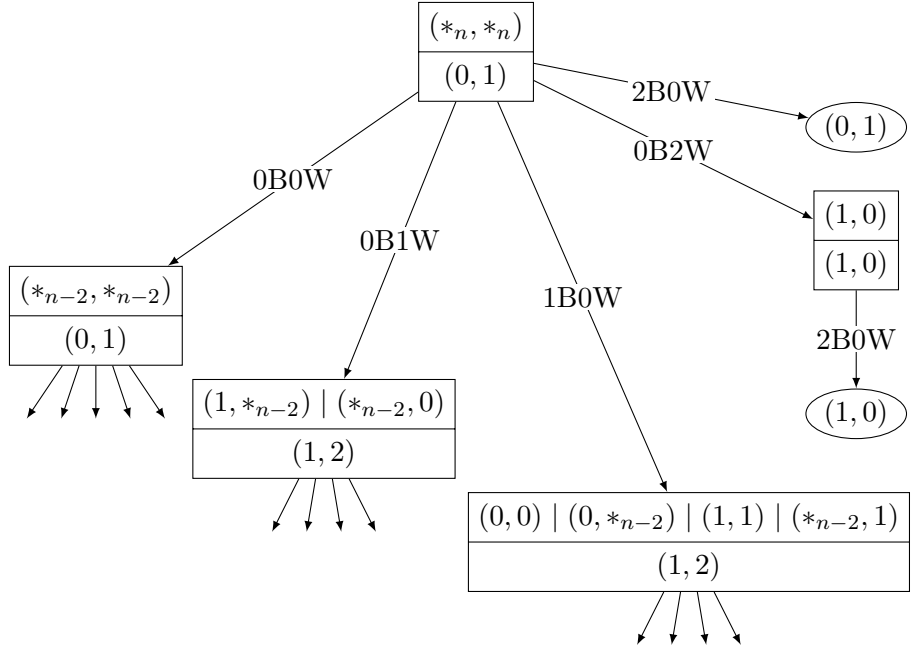

A beginning fragment of the Mastermind game tree for two pegs is shown in
Figure \ref{fig:three}.
An internal node is drawn in rectangle and it is labelled with a maset pattern
and a question.
A leaf is drawn in ellipse and it is labelled with the guessed secret.
The game begins with the full maset pattern $(*_n,*_n)$.
The codebreaker asks the question $(0,1)$.
\begin{itemize}
\item If the codemaker answers with $\ans{0}{0}$, then the codebreaker knows
that the secret does not contain colors 0 and 1.
After normalization, the game state is represented by the maset pattern
$(*_{n-2},*_{n-2})$, which is isomorphic to the initial state, and the
codebreaker asks recursively the question $(0,1)$.
\item If the codemaker answers with $\ans{0}{1}$, then either there is color 1
in position 0, but no colors 0 and 1 in position 1, or color 0 in position 1,
but no colors 0 and 1 in position 0.
The game state is represented by the maset pattern
$(1,*_{n-2}) \mid (*_{n-2},0)$.
The codebreaker asks the question $(1,2)$.
\item If the codemaker answers with $\ans{1}{0}$, then either there is color 0
in position 0, but no color 1 in position 1, or color 1 in position 1, but no
color 0 in position 0.
The game state is represented by the maset pattern
$(0,0) \mid (0,*_{n-2}) \mid (1,1) \mid (*_{n-2},1)$.
The codebreaker asks the question $(1,2)$.
\item If the codemaker answers with $\ans{0}{2}$, then the only possible secret
is $(1,0)$.
The codebreaker asks the question $(1,0)$.
The codemaker must answer with $\ans{2}{0}$, which ends the game.
\item If the codemaker answers with $\ans{2}{0}$, then the secret is guessed and
the game is finished.
\end{itemize}

The number of questions required to solve a maset pattern in the expected case
is $L/N$, where $L$ is the external path length of the game tree for the maset
pattern, i.e.\ the sum of path lengths from the root to each leaf (counted as
the number of edges), and $N$ is the number of secrets in this maset pattern,
i.e.\ the number of leafs in the game tree.
Recursively, external path length of the game tree is the sum of $N$ and the
external path lengths of all subtrees.
The external path length of a leaf is zero.
The external path length of the Mastermind game is the external path length of
the game three with the maset pattern $(*_n,*_n,\dots,*_n)$ in the root.
We obtain the minimal number of question in the expected case when we
minimize the external path length of the game tree.

There are $P=p(p+3)/2$ possible answers in the game with $p$ pegs.
We number answers from 0 to $P-1$, namely 0: $\ans{0}{0}$, 1: $\ans{0}{1}$,
2: $\ans{1}{0}$, 3: $\ans{0}{2}$, 4: $\ans{1}{1}$, 5: $\ans{2}{0}$, \dots,
$P-3$: $\ans{(p-3)}{3}$, $P-2$: $\ans{(p-2)}{2}$, $P-1$: $\ans{p}{0}$.
Note, that the answer $\ans{(p-1)}{1}$ is impossible.

\section{Algorithm}\label{sec:algo}

In this section, we present an algorithm which generates: (a) maset patterns for
Mastermind or AB game, (b) recursive equations for counting the external path
length of the Mastermind or AB game tree.
The algorithm is parametrized with the number of pegs~$p$.

\paragraph{Main Procedure.}
The procedure uses a maset pattern queue.
Initially, the queue contains the maset pattern with $p$ stars,
i.e.\ $M_{p,0}=(*_n,*_n,*_n,\dots,*_n)$.
Let $i$ be the index of the maset pattern in the queue.
Let $s$ denote the actual size of the queue.
The algorithm works as follow.
\begin{itemize}
  \item Initially, we have $i=0$ and $s=1$.
  \item While $i<s$
  \begin{itemize}
    \item we generate all pairwise non-isomorphic questions (denoted by $Q_1$,
      $Q_2$, $Q_3$, \dots) respect to the maset pattern $M_{p,i}$;
    \item for each question $Q_j$
    \begin{itemize}
      \item we split $M_{p,i}$ into maset patterns $M_{p,i,0}$, $M_{p,i,1}$,
        $M_{p,i,2}$, \dots, $M_{p,i,P-1}$, see \textbf{Subprocedure 1};
      \item for $k=0,1,2,\dots,P-1$
      \begin{itemize}
        \item if $M_{p,i,k}$ contains a clause with a star, we normalize
          $M_{p,i,k}$, and then, if the queue does not contain a maset pattern
          isomorphic to $M_{p,i,k}$, we add $M_{p,i,k}$ at the end of the queue
          as $M_{p,s}$ and we increase the size of the queue: $s := s + 1$;
        \item else ($M_{p,i,k}$ does not contain a clause with a star) we
          compute the minimal external path length of the game tree for the
          game starting from this maset pattern, see \textbf{Subprocedure 2};
      \end{itemize}
      \item we compose a recursive equation for the external path length of the
        game tree, when the game is started from the maset pattern $M_{p,i}$ and
        the question $Q_j$ is asked, see \textbf{Subprocedure 3};
        we assign the number $(g.p.i.j)$ to the equation, where $g$ is the game
        name abbreviation: MM or AB;
    \end{itemize}
    \item we increase the maset pattern index in the queue: $i := i + 1$.
  \end{itemize}
\end{itemize}

\paragraph{Subprocedure 1.}
The subprocedure takes the maset pattern $M_{p,i}$ and the question $Q_j$.
It returns the sequence of the maset patterns $M_{p,i,k}$ for
$k=0,1,2,\dots,P-1$, where $M_{p,i,k}$ is the maset pattern, when the
codebreaker asks the question $Q_j$ in the game state $M_{p,i}$ and the
codemaker gives the answer number $k$.
\begin{itemize}
  \item Let $u$ be the number of explicit colors in $M_{p,i}$.
    Let $v$ be such that $v-1$ is the maximal color in $Q_j$.
  \item If $v > u$, we extend $M_{p,i}$ to $v$ colors.
  \item For each $k=0,1,2,\dots,P-1$, let $M_{p,i,k}$ be the empty maset pattern.
  \item For each clause in $M_{p,i}$
  \begin{itemize}
    \item we compute an answer to $Q_j$ using the rules:
      explicit colors are treated as usually,
      any explicit color does not match the star,
      any explicit color does not match the additional color,
      the additional color does not match the star;
    \item if the answer number is $k$, we add the clause to $M_{p,i,k}$.
  \end{itemize}
\end{itemize}

\paragraph{Subprocedure 2.}
The subprocedure takes the maset pattern $M$, which does not contain a star in
clauses.
It returns the minimal external path length $L$ of the game tree for the game
starting from the maset pattern $M$.
This subprocedure is called with $M=M_{p,i,k}$.
\begin{itemize}
  \item If $M$ is the result of the answer number $P-1$, i.e.\ the answer
    $\ans{p}{0}$, then $L=0$;
  \item else if $M$ is empty, then $L=0$;
  \item else if $M$ contains one clause (secret), then $L=1$;
  \item else if $M$ contains two clauses (secrets), then $L=3$
    (see Lemma 1 of \cite{ChL04mm} or Lemma 2 of~\cite{ChL04ab});
  \item else
  \begin{itemize}
    \item we generate all pairwise non-isomorphic questions (denoted by $Q_1$,
      $Q_2$, $Q_3$, \dots) respect to $M$;
    \item for each question $Q_j$, we split $M$ into $M_{j,0}$, $M_{j,1}$,
      $M_{j,2}$, \dots, $M_{j,P-1}$, where $M_{j,k}$ is the maset pattern, when
      the codebreaker asks the question $Q_j$ in the game state $M$ and the
      codemaker gives the answer number $k$;
    \item for each $M_{j,k}$, we compute the minimal external path length
      $L_{j,k}$ using this algorithm recursively;
    \item for each question $Q_j$, we compute the external path length
      $L_j = |M| + L_{j,0} + L_{j,1} + L_{j,2} + \dots + L_{j,P-1}$, where $|M|$
      is the size of the maset $M$, i.e.\ the number of secrets, which is here
      equal to the number of clauses;
    \item we compute the minimum $L = \min_j\{L_j\}$.
  \end{itemize}
\end{itemize}

\paragraph{Subprocedure 3.}
The subprocedure takes the maset pattern $M_{p,i}$ and the sequence of the maset
patterns $M_{p,i,k}$ for $k=0,1,2,\dots,P-1$, obtained after the question $Q_j$.
It returns a recursive equation for the external path length $A_{p,i}(n)$ of the
game tree, when the game is started from the maset pattern $M_{p,i}$, played
with $n$ colors in secrets, the additional color in questions, and when the
question $Q_j$ is asked.
Note that each question $Q_j$ gives a separate equation for $A_{p,i}(n)$.
Note also that $A_{p,0}(n)$ is the external path length of the Mastermind or AB
game tree.

Left side of the equation is the symbol $A_{p,i}(n)$ of the function counting
the external path length.
Right side of the equation is a sum of terms.
If $M_{p,i,k}$ contains the star $*_{n-t}$ and $u$ explicit colors and it is
isomorphic to $l$th maset pattern in the queue, i.e.\ $M_{p,l}$, then we add the
term $A_{p,l}(n-r)$, where $r=t-u$.
The last term is a polynomial $W(n-t)$, being the sum of the number of secrets
represented by $M_{p,i}$ and the sum of the minimal external path lengths for
all $M_{p,i,k}$ which do not contain a star, computed earlier by
\textbf{Subprocedure 2}.

\paragraph{Solving equations.}
For each game state $M_{p,i}$, we deduce, which equation gives the minimal value
of $A_{p,i}(n)$, i.e.\ which question is optimal in this game state.
It could be that there is a question optimal for every $n$, but it could also be
that selection of the optimal question depends on~$n$.

\bigskip
\noindent
The algorithm is implemented in C++.
A compressed archive with the complete source code of the program can be
downloaded from \cite{AVEcode}.
We use the \emph{nauty} package \cite{McK98} to recognize isomorphism.
Equations are solved manually.

\section{Mastermind}\label{sec:mm}

In this section, we present results of computations for Mastermind.
Next, we analyze the results.
The program outputs \LaTeX\ listings of the following form:

\begin{itemize}
  \item[] $M_{p,i} = \mbox{maset pattern}$
    \begin{itemize}
      \item[] question
        \begin{itemize}
          \item[] $\mbox{answer} \rightarrow \mbox{maset pattern}\ [\Rightarrow \mbox{maset pattern}]\ [\Rrightarrow \mbox{maset pattern}]$
          \item[] $\mbox{answer} \rightarrow \mbox{maset pattern}\ [\Rightarrow \mbox{maset pattern}]\ [\Rrightarrow \mbox{maset pattern}]$
            \begin{equation*} \mbox{equation for the external path length} \end{equation*}
        \end{itemize}
      \item[] question
        \begin{itemize}
          \item[] $\mbox{answer} \rightarrow \mbox{maset pattern}\ [\Rightarrow \mbox{maset pattern}]\ [\Rrightarrow \mbox{maset pattern}]$
          \item[] $\mbox{answer} \rightarrow \mbox{maset pattern}\ [\Rightarrow \mbox{maset pattern}]\ [\Rrightarrow \mbox{maset pattern}]$
            \begin{equation*} \mbox{equation for the external path length} \end{equation*}
        \end{itemize}
    \end{itemize}
\end{itemize}

There are printed all found maset patterns pushed in the queue.
For each maset pattern, there is printed a list of pairwise non-isomorphic
questions.
For each question, there is printed a list of all answers.
Sign ``$\rightarrow$'' shows the maset pattern after the answer.
If the answer is impossible, then there is printed the empty maset pattern
$\emptyset$.
Sign ``$\Rightarrow$'' shows the maset pattern after normalization.
It is printed only if normalization changes the maset pattern.
Sign ``$\Rrightarrow$'' shows the isomorphic maset pattern in the queue.
It is printed only if the maset pattern in the queue is different than the maset
pattern after normalization.
For each question, there is printed a recursive equation for the external path
length of the game tree.

Next subsections contain results for one peg, two pegs and three pegs.
The result for one peg is quite trivial, but it is helpful to check correctness
of the program.

\subsection{One peg}

\begin{itemize}
  \item[] $M_{1,0} = (*_{n})$
    \begin{itemize}
      \item[] $(0)$
        \begin{itemize}
          \item[] $\ans{0}{0} \rightarrow (*_{n-1})$
          \item[] $\ans{1}{0} \rightarrow (0)$
            \begin{equation} A_{1,0}(n) = A_{1,0}(n-1) + n \tag{MM.1.0.1} \label{eq:mm.1.0.1} \end{equation}
        \end{itemize}
    \end{itemize}
\end{itemize}

% 1 maset-pattern
% 1 equation
% dependency cycle not found

One maset pattern is found: $(*_n)$.
It represents $n$ secrets.
Hence, we have that $A_{1,0}(1) = 1$.
Therefore, the solution of equation (\ref{eq:mm.1.0.1}) is
$A_{1,0}(n)=\frac{1}{2}(n+1)n$, and the minimal expected number of questions for
one peg Mastermind is $A_{1,0}(n)/n=\frac{1}{2}(n+1)$, which is consistent with
intuition.

\subsection{Two pegs}

\begin{itemize}
  \item[] $M_{2,0} = (*_{n},*_{n})$
    \begin{itemize}
      \item[] $(0,0)$
        \begin{itemize}
          \item[] $\ans{0}{0} \rightarrow (*_{n-1},*_{n-1})$
          \item[] $\ans{0}{1} \rightarrow \emptyset$
          \item[] $\ans{1}{0} \rightarrow (0,*_{n-1}) \mid (*_{n-1},0)$
          \item[] $\ans{0}{2} \rightarrow \emptyset$
          \item[] $\ans{2}{0} \rightarrow (0,0)$
            \begin{equation} A_{2,0}(n) = A_{2,0}(n-1) + A_{2,1}(n) + n^2 \tag{MM.2.0.1} \label{eq:mm.2.0.1} \end{equation}
        \end{itemize}
      \item[] $(0,1)$
        \begin{itemize}
          \item[] $\ans{0}{0} \rightarrow (*_{n-2},*_{n-2})$
          \item[] $\ans{0}{1} \rightarrow (1,*_{n-2}) \mid (*_{n-2},0)$
          \item[] $\ans{1}{0} \rightarrow (0,0) \mid (0,*_{n-2}) \mid (1,1) \mid (*_{n-2},1)$
          \item[] $\ans{0}{2} \rightarrow (1,0)$
          \item[] $\ans{2}{0} \rightarrow (0,1)$
            \begin{equation} A_{2,0}(n) = A_{2,0}(n-2) + A_{2,2}(n) + A_{2,3}(n) + n^2 + 1 \tag{MM.2.0.2} \label{eq:mm.2.0.2} \end{equation}
        \end{itemize}
      \item[] $(0,\ac)$
        \begin{itemize}
          \item[] $\ans{0}{0} \rightarrow (*_{n-1},*_{n-1})$
          \item[] $\ans{0}{1} \rightarrow (*_{n-1},0)$
          \item[] $\ans{1}{0} \rightarrow (0,0) \mid (0,*_{n-1})$
          \item[] $\ans{0}{2} \rightarrow \emptyset$
          \item[] $\ans{2}{0} \rightarrow \emptyset$
            \begin{equation} A_{2,0}(n) = A_{2,0}(n-1) + A_{2,4}(n) + A_{2,5}(n) + n^2 \tag{MM.2.0.3} \label{eq:mm.2.0.3} \end{equation}
        \end{itemize}
    \end{itemize}
  \item[] $M_{2,1} = (0,*_{n-1}) \mid (*_{n-1},0)$
    \begin{itemize}
      \item[] $(0,0)$
        \begin{itemize}
          \item[] $\ans{0}{0} \rightarrow \emptyset$
          \item[] $\ans{0}{1} \rightarrow \emptyset$
          \item[] $\ans{1}{0} \rightarrow (0,*_{n-1}) \mid (*_{n-1},0)$
          \item[] $\ans{0}{2} \rightarrow \emptyset$
          \item[] $\ans{2}{0} \rightarrow \emptyset$
            \begin{equation} A_{2,1}(n) = A_{2,1}(n) + 2(n-1) \tag{MM.2.1.1} \label{eq:mm.2.1.1} \end{equation}
        \end{itemize}
      \item[] $(0,1)$
        \begin{itemize}
          \item[] $\ans{0}{0} \rightarrow \emptyset$
          \item[] $\ans{0}{1} \rightarrow (*_{n-2},0)$
          \item[] $\ans{1}{0} \rightarrow (0,*_{n-2}) \Rrightarrow (*_{n-2},0)$
          \item[] $\ans{0}{2} \rightarrow (1,0)$
          \item[] $\ans{2}{0} \rightarrow (0,1)$
            \begin{equation} A_{2,1}(n) = A_{2,4}(n-1) + A_{2,4}(n-1) + 2(n-1) + 1 \tag{MM.2.1.2} \label{eq:mm.2.1.2} \end{equation}
        \end{itemize}
      \item[] $(0,\ac)$
        \begin{itemize}
          \item[] $\ans{0}{0} \rightarrow \emptyset$
          \item[] $\ans{0}{1} \rightarrow (*_{n-1},0)$
          \item[] $\ans{1}{0} \rightarrow (0,*_{n-1}) \Rrightarrow (*_{n-1},0)$
          \item[] $\ans{0}{2} \rightarrow \emptyset$
          \item[] $\ans{2}{0} \rightarrow \emptyset$
            \begin{equation} A_{2,1}(n) = A_{2,4}(n) + A_{2,4}(n) + 2(n-1) \tag{MM.2.1.3} \label{eq:mm.2.1.3} \end{equation}
        \end{itemize}
      \item[] $(1,1)$
        \begin{itemize}
          \item[] $\ans{0}{0} \rightarrow (0,*_{n-2}) \mid (*_{n-2},0)$
          \item[] $\ans{0}{1} \rightarrow \emptyset$
          \item[] $\ans{1}{0} \rightarrow (0,1) \mid (1,0)$
          \item[] $\ans{0}{2} \rightarrow \emptyset$
          \item[] $\ans{2}{0} \rightarrow \emptyset$
            \begin{equation} A_{2,1}(n) = A_{2,1}(n-1) + 2(n-1) + 3 \tag{MM.2.1.4} \label{eq:mm.2.1.4} \end{equation}
        \end{itemize}
      \item[] $(1,2)$
        \begin{itemize}
          \item[] $\ans{0}{0} \rightarrow (0,*_{n-3}) \mid (*_{n-3},0)$
          \item[] $\ans{0}{1} \rightarrow (0,1) \mid (2,0)$
          \item[] $\ans{1}{0} \rightarrow (0,2) \mid (1,0)$
          \item[] $\ans{0}{2} \rightarrow \emptyset$
          \item[] $\ans{2}{0} \rightarrow \emptyset$
            \begin{equation} A_{2,1}(n) = A_{2,1}(n-2) + 2(n-1) + 6 \tag{MM.2.1.5} \label{eq:mm.2.1.5} \end{equation}
        \end{itemize}
      \item[] $(1,\ac)$
        \begin{itemize}
          \item[] $\ans{0}{0} \rightarrow (0,*_{n-2}) \mid (*_{n-2},0)$
          \item[] $\ans{0}{1} \rightarrow (0,1)$
          \item[] $\ans{1}{0} \rightarrow (1,0)$
          \item[] $\ans{0}{2} \rightarrow \emptyset$
          \item[] $\ans{2}{0} \rightarrow \emptyset$
            \begin{equation} A_{2,1}(n) = A_{2,1}(n-1) + 2(n-1) + 2 \tag{MM.2.1.6} \label{eq:mm.2.1.6} \end{equation}
        \end{itemize}
    \end{itemize}
  \item[] $M_{2,2} = (1,*_{n-2}) \mid (*_{n-2},0)$
    \begin{itemize}
      \item[] $(0,0)$
        \begin{itemize}
          \item[] $\ans{0}{0} \rightarrow (1,*_{n-2}) \Rightarrow (0,*_{n-2}) \Rrightarrow (*_{n-2},0)$
          \item[] $\ans{0}{1} \rightarrow \emptyset$
          \item[] $\ans{1}{0} \rightarrow (*_{n-2},0)$
          \item[] $\ans{0}{2} \rightarrow \emptyset$
          \item[] $\ans{2}{0} \rightarrow \emptyset$
            \begin{equation} A_{2,2}(n) = A_{2,4}(n-1) + A_{2,4}(n-1) + 2(n-2) \tag{MM.2.2.1} \label{eq:mm.2.2.1} \end{equation}
        \end{itemize}
      \item[] $(0,1)$
        \begin{itemize}
          \item[] $\ans{0}{0} \rightarrow \emptyset$
          \item[] $\ans{0}{1} \rightarrow (1,*_{n-2}) \mid (*_{n-2},0)$
          \item[] $\ans{1}{0} \rightarrow \emptyset$
          \item[] $\ans{0}{2} \rightarrow \emptyset$
          \item[] $\ans{2}{0} \rightarrow \emptyset$
            \begin{equation} A_{2,2}(n) = A_{2,2}(n) + 2(n-2) \tag{MM.2.2.2} \label{eq:mm.2.2.2} \end{equation}
        \end{itemize}
      \item[] $(0,2)$
        \begin{itemize}
          \item[] $\ans{0}{0} \rightarrow (1,*_{n-3}) \Rightarrow (0,*_{n-3}) \Rrightarrow (*_{n-3},0)$
          \item[] $\ans{0}{1} \rightarrow (*_{n-3},0)$
          \item[] $\ans{1}{0} \rightarrow (1,2)$
          \item[] $\ans{0}{2} \rightarrow (2,0)$
          \item[] $\ans{2}{0} \rightarrow \emptyset$
            \begin{equation} A_{2,2}(n) = A_{2,4}(n-2) + A_{2,4}(n-2) + 2(n-2) + 2 \tag{MM.2.2.3} \label{eq:mm.2.2.3} \end{equation}
        \end{itemize}
      \item[] $(0,\ac)$
        \begin{itemize}
          \item[] $\ans{0}{0} \rightarrow (1,*_{n-2}) \Rightarrow (0,*_{n-2}) \Rrightarrow (*_{n-2},0)$
          \item[] $\ans{0}{1} \rightarrow (*_{n-2},0)$
          \item[] $\ans{1}{0} \rightarrow \emptyset$
          \item[] $\ans{0}{2} \rightarrow \emptyset$
          \item[] $\ans{2}{0} \rightarrow \emptyset$
            \begin{equation} A_{2,2}(n) = A_{2,4}(n-1) + A_{2,4}(n-1) + 2(n-2) \tag{MM.2.2.4} \label{eq:mm.2.2.4} \end{equation}
        \end{itemize}
      \item[] $(1,0)$
        \begin{itemize}
          \item[] $\ans{0}{0} \rightarrow \emptyset$
          \item[] $\ans{0}{1} \rightarrow \emptyset$
          \item[] $\ans{1}{0} \rightarrow (1,*_{n-2}) \mid (*_{n-2},0)$
          \item[] $\ans{0}{2} \rightarrow \emptyset$
          \item[] $\ans{2}{0} \rightarrow \emptyset$
            \begin{equation} A_{2,2}(n) = A_{2,2}(n) + 2(n-2) \tag{MM.2.2.5} \label{eq:mm.2.2.5} \end{equation}
        \end{itemize}
      \item[] $(1,2)$
        \begin{itemize}
          \item[] $\ans{0}{0} \rightarrow (*_{n-3},0)$
          \item[] $\ans{0}{1} \rightarrow (2,0)$
          \item[] $\ans{1}{0} \rightarrow (1,*_{n-3}) \Rightarrow (0,*_{n-3}) \Rrightarrow (*_{n-3},0)$
          \item[] $\ans{0}{2} \rightarrow \emptyset$
          \item[] $\ans{2}{0} \rightarrow (1,2)$
            \begin{equation} A_{2,2}(n) = A_{2,4}(n-2) + A_{2,4}(n-2) + 2(n-2) + 1 \tag{MM.2.2.6} \label{eq:mm.2.2.6} \end{equation}
        \end{itemize}
      \item[] $(1,\ac)$
        \begin{itemize}
          \item[] $\ans{0}{0} \rightarrow (*_{n-2},0)$
          \item[] $\ans{0}{1} \rightarrow \emptyset$
          \item[] $\ans{1}{0} \rightarrow (1,*_{n-2}) \Rightarrow (0,*_{n-2}) \Rrightarrow (*_{n-2},0)$
          \item[] $\ans{0}{2} \rightarrow \emptyset$
          \item[] $\ans{2}{0} \rightarrow \emptyset$
            \begin{equation} A_{2,2}(n) = A_{2,4}(n-1) + A_{2,4}(n-1) + 2(n-2) \tag{MM.2.2.7} \label{eq:mm.2.2.7} \end{equation}
        \end{itemize}
      \item[] $(2,2)$
        \begin{itemize}
          \item[] $\ans{0}{0} \rightarrow (1,*_{n-3}) \mid (*_{n-3},0)$
          \item[] $\ans{0}{1} \rightarrow \emptyset$
          \item[] $\ans{1}{0} \rightarrow (1,2) \mid (2,0)$
          \item[] $\ans{0}{2} \rightarrow \emptyset$
          \item[] $\ans{2}{0} \rightarrow \emptyset$
            \begin{equation} A_{2,2}(n) = A_{2,2}(n-1) + 2(n-2) + 3 \tag{MM.2.2.8} \label{eq:mm.2.2.8} \end{equation}
        \end{itemize}
      \item[] $(2,3)$
        \begin{itemize}
          \item[] $\ans{0}{0} \rightarrow (1,*_{n-4}) \mid (*_{n-4},0)$
          \item[] $\ans{0}{1} \rightarrow (1,2) \mid (3,0)$
          \item[] $\ans{1}{0} \rightarrow (1,3) \mid (2,0)$
          \item[] $\ans{0}{2} \rightarrow \emptyset$
          \item[] $\ans{2}{0} \rightarrow \emptyset$
            \begin{equation} A_{2,2}(n) = A_{2,2}(n-2) + 2(n-2) + 6 \tag{MM.2.2.9} \label{eq:mm.2.2.9} \end{equation}
        \end{itemize}
      \item[] $(2,\ac)$
        \begin{itemize}
          \item[] $\ans{0}{0} \rightarrow (1,*_{n-3}) \mid (*_{n-3},0)$
          \item[] $\ans{0}{1} \rightarrow (1,2)$
          \item[] $\ans{1}{0} \rightarrow (2,0)$
          \item[] $\ans{0}{2} \rightarrow \emptyset$
          \item[] $\ans{2}{0} \rightarrow \emptyset$
            \begin{equation} A_{2,2}(n) = A_{2,2}(n-1) + 2(n-2) + 2 \tag{MM.2.2.10} \label{eq:mm.2.2.10} \end{equation}
        \end{itemize}
    \end{itemize}
  \item[] $M_{2,3} = (0,0) \mid (0,*_{n-2}) \mid (1,1) \mid (*_{n-2},1)$
    \begin{itemize}
      \item[] $(0,0)$
        \begin{itemize}
          \item[] $\ans{0}{0} \rightarrow (1,1) \mid (*_{n-2},1) \Rightarrow (0,0) \mid (*_{n-2},0) \Rrightarrow (0,0) \mid (0,*_{n-2})$
          \item[] $\ans{0}{1} \rightarrow \emptyset$
          \item[] $\ans{1}{0} \rightarrow (0,*_{n-2}) \Rrightarrow (*_{n-2},0)$
          \item[] $\ans{0}{2} \rightarrow \emptyset$
          \item[] $\ans{2}{0} \rightarrow (0,0)$
            \begin{equation} A_{2,3}(n) = A_{2,5}(n-1) + A_{2,4}(n-1) + 2(n-2) + 2 \tag{MM.2.3.1} \label{eq:mm.2.3.1} \end{equation}
        \end{itemize}
      \item[] $(0,1)$
        \begin{itemize}
          \item[] $\ans{0}{0} \rightarrow \emptyset$
          \item[] $\ans{0}{1} \rightarrow \emptyset$
          \item[] $\ans{1}{0} \rightarrow (0,0) \mid (0,*_{n-2}) \mid (1,1) \mid (*_{n-2},1)$
          \item[] $\ans{0}{2} \rightarrow \emptyset$
          \item[] $\ans{2}{0} \rightarrow \emptyset$
            \begin{equation} A_{2,3}(n) = A_{2,3}(n) + 2(n-2) + 2 \tag{MM.2.3.2} \label{eq:mm.2.3.2} \end{equation}
        \end{itemize}
      \item[] $(0,2)$
        \begin{itemize}
          \item[] $\ans{0}{0} \rightarrow (1,1) \mid (*_{n-3},1) \Rightarrow (0,0) \mid (*_{n-3},0) \Rrightarrow (0,0) \mid (0,*_{n-3})$
          \item[] $\ans{0}{1} \rightarrow (2,1)$
          \item[] $\ans{1}{0} \rightarrow (0,0) \mid (0,*_{n-3})$
          \item[] $\ans{0}{2} \rightarrow \emptyset$
          \item[] $\ans{2}{0} \rightarrow (0,2)$
            \begin{equation} A_{2,3}(n) = A_{2,5}(n-2) + A_{2,5}(n-2) + 2(n-2) + 3 \tag{MM.2.3.3} \label{eq:mm.2.3.3} \end{equation}
        \end{itemize}
      \item[] $(0,\ac)$
        \begin{itemize}
          \item[] $\ans{0}{0} \rightarrow (1,1) \mid (*_{n-2},1) \Rightarrow (0,0) \mid (*_{n-2},0) \Rrightarrow (0,0) \mid (0,*_{n-2})$
          \item[] $\ans{0}{1} \rightarrow \emptyset$
          \item[] $\ans{1}{0} \rightarrow (0,0) \mid (0,*_{n-2})$
          \item[] $\ans{0}{2} \rightarrow \emptyset$
          \item[] $\ans{2}{0} \rightarrow \emptyset$
            \begin{equation} A_{2,3}(n) = A_{2,5}(n-1) + A_{2,5}(n-1) + 2(n-2) + 2 \tag{MM.2.3.4} \label{eq:mm.2.3.4} \end{equation}
        \end{itemize}
      \item[] $(1,0)$
        \begin{itemize}
          \item[] $\ans{0}{0} \rightarrow \emptyset$
          \item[] $\ans{0}{1} \rightarrow (0,*_{n-2}) \mid (*_{n-2},1) \Rrightarrow (1,*_{n-2}) \mid (*_{n-2},0)$
          \item[] $\ans{1}{0} \rightarrow (0,0) \mid (1,1)$
          \item[] $\ans{0}{2} \rightarrow \emptyset$
          \item[] $\ans{2}{0} \rightarrow \emptyset$
            \begin{equation} A_{2,3}(n) = A_{2,2}(n) + 2(n-2) + 5 \tag{MM.2.3.5} \label{eq:mm.2.3.5} \end{equation}
        \end{itemize}
      \item[] $(1,2)$
        \begin{itemize}
          \item[] $\ans{0}{0} \rightarrow (0,0) \mid (0,*_{n-3})$
          \item[] $\ans{0}{1} \rightarrow (*_{n-3},1) \Rightarrow (*_{n-3},0)$
          \item[] $\ans{1}{0} \rightarrow (0,2) \mid (1,1)$
          \item[] $\ans{0}{2} \rightarrow (2,1)$
          \item[] $\ans{2}{0} \rightarrow \emptyset$
            \begin{equation} A_{2,3}(n) = A_{2,5}(n-2) + A_{2,4}(n-2) + 2(n-2) + 6 \tag{MM.2.3.6} \label{eq:mm.2.3.6} \end{equation}
        \end{itemize}
      \item[] $(1,\ac)$
        \begin{itemize}
          \item[] $\ans{0}{0} \rightarrow (0,0) \mid (0,*_{n-2})$
          \item[] $\ans{0}{1} \rightarrow (*_{n-2},1) \Rightarrow (*_{n-2},0)$
          \item[] $\ans{1}{0} \rightarrow (1,1)$
          \item[] $\ans{0}{2} \rightarrow \emptyset$
          \item[] $\ans{2}{0} \rightarrow \emptyset$
            \begin{equation} A_{2,3}(n) = A_{2,5}(n-1) + A_{2,4}(n-1) + 2(n-2) + 3 \tag{MM.2.3.7} \label{eq:mm.2.3.7} \end{equation}
        \end{itemize}
      \item[] $(2,2)$
        \begin{itemize}
          \item[] $\ans{0}{0} \rightarrow (0,0) \mid (0,*_{n-3}) \mid (1,1) \mid (*_{n-3},1)$
          \item[] $\ans{0}{1} \rightarrow \emptyset$
          \item[] $\ans{1}{0} \rightarrow (0,2) \mid (2,1)$
          \item[] $\ans{0}{2} \rightarrow \emptyset$
          \item[] $\ans{2}{0} \rightarrow \emptyset$
            \begin{equation} A_{2,3}(n) = A_{2,3}(n-1) + 2(n-2) + 5 \tag{MM.2.3.8} \label{eq:mm.2.3.8} \end{equation}
        \end{itemize}
      \item[] $(2,3)$
        \begin{itemize}
          \item[] $\ans{0}{0} \rightarrow (0,0) \mid (0,*_{n-4}) \mid (1,1) \mid (*_{n-4},1)$
          \item[] $\ans{0}{1} \rightarrow (0,2) \mid (3,1)$
          \item[] $\ans{1}{0} \rightarrow (0,3) \mid (2,1)$
          \item[] $\ans{0}{2} \rightarrow \emptyset$
          \item[] $\ans{2}{0} \rightarrow \emptyset$
            \begin{equation} A_{2,3}(n) = A_{2,3}(n-2) + 2(n-2) + 8 \tag{MM.2.3.9} \label{eq:mm.2.3.9} \end{equation}
        \end{itemize}
      \item[] $(2,\ac)$
        \begin{itemize}
          \item[] $\ans{0}{0} \rightarrow (0,0) \mid (0,*_{n-3}) \mid (1,1) \mid (*_{n-3},1)$
          \item[] $\ans{0}{1} \rightarrow (0,2)$
          \item[] $\ans{1}{0} \rightarrow (2,1)$
          \item[] $\ans{0}{2} \rightarrow \emptyset$
          \item[] $\ans{2}{0} \rightarrow \emptyset$
            \begin{equation} A_{2,3}(n) = A_{2,3}(n-1) + 2(n-2) + 4 \tag{MM.2.3.10} \label{eq:mm.2.3.10} \end{equation}
        \end{itemize}
    \end{itemize}
  \item[] $M_{2,4} = (*_{n-1},0)$
    \begin{itemize}
      \item[] $(0,0)$
        \begin{itemize}
          \item[] $\ans{0}{0} \rightarrow \emptyset$
          \item[] $\ans{0}{1} \rightarrow \emptyset$
          \item[] $\ans{1}{0} \rightarrow (*_{n-1},0)$
          \item[] $\ans{0}{2} \rightarrow \emptyset$
          \item[] $\ans{2}{0} \rightarrow \emptyset$
            \begin{equation} A_{2,4}(n) = A_{2,4}(n) + (n-1) \tag{MM.2.4.1} \label{eq:mm.2.4.1} \end{equation}
        \end{itemize}
      \item[] $(0,1)$
        \begin{itemize}
          \item[] $\ans{0}{0} \rightarrow \emptyset$
          \item[] $\ans{0}{1} \rightarrow (*_{n-2},0)$
          \item[] $\ans{1}{0} \rightarrow \emptyset$
          \item[] $\ans{0}{2} \rightarrow (1,0)$
          \item[] $\ans{2}{0} \rightarrow \emptyset$
            \begin{equation} A_{2,4}(n) = A_{2,4}(n-1) + (n-1) + 1 \tag{MM.2.4.2} \label{eq:mm.2.4.2} \end{equation}
        \end{itemize}
      \item[] $(0,\ac)$
        \begin{itemize}
          \item[] $\ans{0}{0} \rightarrow \emptyset$
          \item[] $\ans{0}{1} \rightarrow (*_{n-1},0)$
          \item[] $\ans{1}{0} \rightarrow \emptyset$
          \item[] $\ans{0}{2} \rightarrow \emptyset$
          \item[] $\ans{2}{0} \rightarrow \emptyset$
            \begin{equation} A_{2,4}(n) = A_{2,4}(n) + (n-1) \tag{MM.2.4.3} \label{eq:mm.2.4.3} \end{equation}
        \end{itemize}
      \item[] $(1,0)$
        \begin{itemize}
          \item[] $\ans{0}{0} \rightarrow \emptyset$
          \item[] $\ans{0}{1} \rightarrow \emptyset$
          \item[] $\ans{1}{0} \rightarrow (*_{n-2},0)$
          \item[] $\ans{0}{2} \rightarrow \emptyset$
          \item[] $\ans{2}{0} \rightarrow (1,0)$
            \begin{equation} A_{2,4}(n) = A_{2,4}(n-1) + (n-1) \tag{MM.2.4.4} \label{eq:mm.2.4.4} \end{equation}
        \end{itemize}
      \item[] $(1,1)$
        \begin{itemize}
          \item[] $\ans{0}{0} \rightarrow (*_{n-2},0)$
          \item[] $\ans{0}{1} \rightarrow \emptyset$
          \item[] $\ans{1}{0} \rightarrow (1,0)$
          \item[] $\ans{0}{2} \rightarrow \emptyset$
          \item[] $\ans{2}{0} \rightarrow \emptyset$
            \begin{equation} A_{2,4}(n) = A_{2,4}(n-1) + (n-1) + 1 \tag{MM.2.4.5} \label{eq:mm.2.4.5} \end{equation}
        \end{itemize}
      \item[] $(1,2)$
        \begin{itemize}
          \item[] $\ans{0}{0} \rightarrow (*_{n-3},0)$
          \item[] $\ans{0}{1} \rightarrow (2,0)$
          \item[] $\ans{1}{0} \rightarrow (1,0)$
          \item[] $\ans{0}{2} \rightarrow \emptyset$
          \item[] $\ans{2}{0} \rightarrow \emptyset$
            \begin{equation} A_{2,4}(n) = A_{2,4}(n-2) + (n-1) + 2 \tag{MM.2.4.6} \label{eq:mm.2.4.6} \end{equation}
        \end{itemize}
      \item[] $(1,\ac)$
        \begin{itemize}
          \item[] $\ans{0}{0} \rightarrow (*_{n-2},0)$
          \item[] $\ans{0}{1} \rightarrow \emptyset$
          \item[] $\ans{1}{0} \rightarrow (1,0)$
          \item[] $\ans{0}{2} \rightarrow \emptyset$
          \item[] $\ans{2}{0} \rightarrow \emptyset$
            \begin{equation} A_{2,4}(n) = A_{2,4}(n-1) + (n-1) + 1 \tag{MM.2.4.7} \label{eq:mm.2.4.7} \end{equation}
        \end{itemize}
      \item[] $(\ac,0)$
        \begin{itemize}
          \item[] $\ans{0}{0} \rightarrow \emptyset$
          \item[] $\ans{0}{1} \rightarrow \emptyset$
          \item[] $\ans{1}{0} \rightarrow (*_{n-1},0)$
          \item[] $\ans{0}{2} \rightarrow \emptyset$
          \item[] $\ans{2}{0} \rightarrow \emptyset$
            \begin{equation} A_{2,4}(n) = A_{2,4}(n) + (n-1) \tag{MM.2.4.8} \label{eq:mm.2.4.8} \end{equation}
        \end{itemize}
      \item[] $(\ac,1)$
        \begin{itemize}
          \item[] $\ans{0}{0} \rightarrow (*_{n-2},0)$
          \item[] $\ans{0}{1} \rightarrow (1,0)$
          \item[] $\ans{1}{0} \rightarrow \emptyset$
          \item[] $\ans{0}{2} \rightarrow \emptyset$
          \item[] $\ans{2}{0} \rightarrow \emptyset$
            \begin{equation} A_{2,4}(n) = A_{2,4}(n-1) + (n-1) + 1 \tag{MM.2.4.9} \label{eq:mm.2.4.9} \end{equation}
        \end{itemize}
    \end{itemize}
  \item[] $M_{2,5} = (0,0) \mid (0,*_{n-1})$
    \begin{itemize}
      \item[] $(0,0)$
        \begin{itemize}
          \item[] $\ans{0}{0} \rightarrow \emptyset$
          \item[] $\ans{0}{1} \rightarrow \emptyset$
          \item[] $\ans{1}{0} \rightarrow (0,*_{n-1}) \Rrightarrow (*_{n-1},0)$
          \item[] $\ans{0}{2} \rightarrow \emptyset$
          \item[] $\ans{2}{0} \rightarrow (0,0)$
            \begin{equation} A_{2,5}(n) = A_{2,4}(n) + (n-1) + 1 \tag{MM.2.5.1} \label{eq:mm.2.5.1} \end{equation}
        \end{itemize}
      \item[] $(0,1)$
        \begin{itemize}
          \item[] $\ans{0}{0} \rightarrow \emptyset$
          \item[] $\ans{0}{1} \rightarrow \emptyset$
          \item[] $\ans{1}{0} \rightarrow (0,0) \mid (0,*_{n-2})$
          \item[] $\ans{0}{2} \rightarrow \emptyset$
          \item[] $\ans{2}{0} \rightarrow (0,1)$
            \begin{equation} A_{2,5}(n) = A_{2,5}(n-1) + (n-1) + 1 \tag{MM.2.5.2} \label{eq:mm.2.5.2} \end{equation}
        \end{itemize}
      \item[] $(0,\ac)$
        \begin{itemize}
          \item[] $\ans{0}{0} \rightarrow \emptyset$
          \item[] $\ans{0}{1} \rightarrow \emptyset$
          \item[] $\ans{1}{0} \rightarrow (0,0) \mid (0,*_{n-1})$
          \item[] $\ans{0}{2} \rightarrow \emptyset$
          \item[] $\ans{2}{0} \rightarrow \emptyset$
            \begin{equation} A_{2,5}(n) = A_{2,5}(n) + (n-1) + 1 \tag{MM.2.5.3} \label{eq:mm.2.5.3} \end{equation}
        \end{itemize}
      \item[] $(1,0)$
        \begin{itemize}
          \item[] $\ans{0}{0} \rightarrow \emptyset$
          \item[] $\ans{0}{1} \rightarrow (0,*_{n-2}) \Rrightarrow (*_{n-2},0)$
          \item[] $\ans{1}{0} \rightarrow (0,0)$
          \item[] $\ans{0}{2} \rightarrow (0,1)$
          \item[] $\ans{2}{0} \rightarrow \emptyset$
            \begin{equation} A_{2,5}(n) = A_{2,4}(n-1) + (n-1) + 3 \tag{MM.2.5.4} \label{eq:mm.2.5.4} \end{equation}
        \end{itemize}
      \item[] $(1,1)$
        \begin{itemize}
          \item[] $\ans{0}{0} \rightarrow (0,0) \mid (0,*_{n-2})$
          \item[] $\ans{0}{1} \rightarrow \emptyset$
          \item[] $\ans{1}{0} \rightarrow (0,1)$
          \item[] $\ans{0}{2} \rightarrow \emptyset$
          \item[] $\ans{2}{0} \rightarrow \emptyset$
            \begin{equation} A_{2,5}(n) = A_{2,5}(n-1) + (n-1) + 2 \tag{MM.2.5.5} \label{eq:mm.2.5.5} \end{equation}
        \end{itemize}
      \item[] $(1,2)$
        \begin{itemize}
          \item[] $\ans{0}{0} \rightarrow (0,0) \mid (0,*_{n-3})$
          \item[] $\ans{0}{1} \rightarrow (0,1)$
          \item[] $\ans{1}{0} \rightarrow (0,2)$
          \item[] $\ans{0}{2} \rightarrow \emptyset$
          \item[] $\ans{2}{0} \rightarrow \emptyset$
            \begin{equation} A_{2,5}(n) = A_{2,5}(n-2) + (n-1) + 3 \tag{MM.2.5.6} \label{eq:mm.2.5.6} \end{equation}
        \end{itemize}
      \item[] $(1,\ac)$
        \begin{itemize}
          \item[] $\ans{0}{0} \rightarrow (0,0) \mid (0,*_{n-2})$
          \item[] $\ans{0}{1} \rightarrow (0,1)$
          \item[] $\ans{1}{0} \rightarrow \emptyset$
          \item[] $\ans{0}{2} \rightarrow \emptyset$
          \item[] $\ans{2}{0} \rightarrow \emptyset$
            \begin{equation} A_{2,5}(n) = A_{2,5}(n-1) + (n-1) + 2 \tag{MM.2.5.7} \label{eq:mm.2.5.7} \end{equation}
        \end{itemize}
      \item[] $(\ac,0)$
        \begin{itemize}
          \item[] $\ans{0}{0} \rightarrow \emptyset$
          \item[] $\ans{0}{1} \rightarrow (0,*_{n-1}) \Rrightarrow (*_{n-1},0)$
          \item[] $\ans{1}{0} \rightarrow (0,0)$
          \item[] $\ans{0}{2} \rightarrow \emptyset$
          \item[] $\ans{2}{0} \rightarrow \emptyset$
            \begin{equation} A_{2,5}(n) = A_{2,4}(n) + (n-1) + 2 \tag{MM.2.5.8} \label{eq:mm.2.5.8} \end{equation}
        \end{itemize}
      \item[] $(\ac,1)$
        \begin{itemize}
          \item[] $\ans{0}{0} \rightarrow (0,0) \mid (0,*_{n-2})$
          \item[] $\ans{0}{1} \rightarrow \emptyset$
          \item[] $\ans{1}{0} \rightarrow (0,1)$
          \item[] $\ans{0}{2} \rightarrow \emptyset$
          \item[] $\ans{2}{0} \rightarrow \emptyset$
            \begin{equation} A_{2,5}(n) = A_{2,5}(n-1) + (n-1) + 2 \tag{MM.2.5.9} \label{eq:mm.2.5.9} \end{equation}
        \end{itemize}
    \end{itemize}
\end{itemize}

% 6 maset-patterns
% 47 equations
% dependency cycle not found

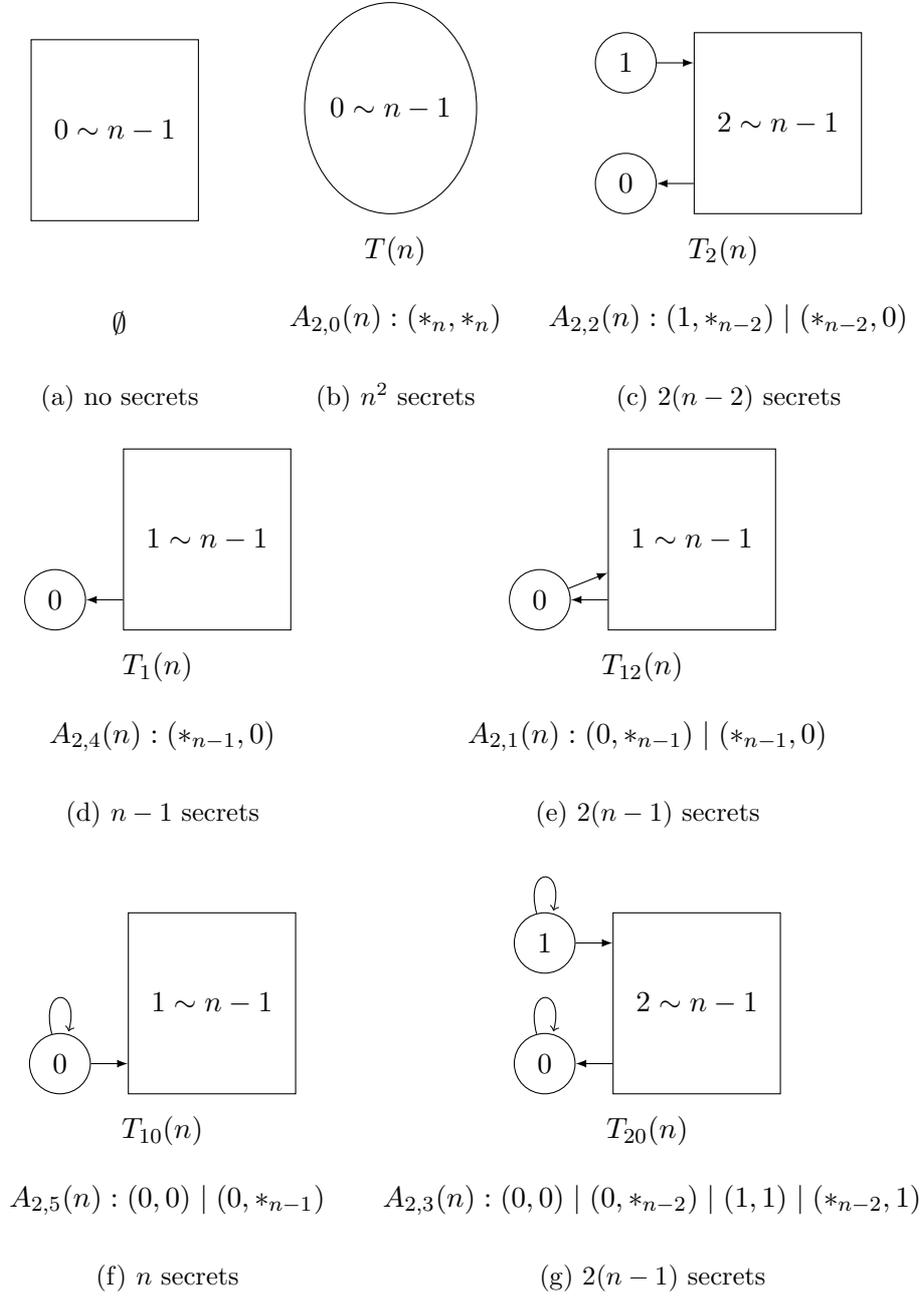
\begin{figure}[htbp]
  \tikzstyle{vertex}=[circle, draw, minimum size=8mm, inner sep=0cm]
  \tikzstyle{emptygraph}=[rectangle, draw, minimum width=22mm, minimum height=24mm, inner sep=0cm]
  \tikzstyle{fullgraph}=[ellipse, draw, minimum width=22mm, minimum height=28mm, inner sep=0cm]
  \begin{center}
    \begin{subfigure}[b]{35mm}
      \begin{center}
        \begin{tikzpicture}
          \node[emptygraph] at (0,0) {$0 \sim n-1$};
        \end{tikzpicture}
        \smallskip\par
        \null\par
        \bigskip\par
        $\emptyset$
      \end{center}
      \caption{no secrets}\label{fig:mm:a}
    \end{subfigure}
    \begin{subfigure}[b]{35mm}
      \begin{center}
        \begin{tikzpicture}
          \node[fullgraph] at (0,0) {$0 \sim n-1$};
        \end{tikzpicture}
        \smallskip\par
        $T(n)$\par
        \bigskip\par
        $A_{2,0}(n): (*_{n},*_{n})$\par
      \end{center}
      \caption{$n^2$ secrets}\label{fig:mm:b}
    \end{subfigure}
    \begin{subfigure}[b]{50mm}
      \begin{center}
        \begin{tikzpicture}
          \node[vertex] (1) at (-2,0.8) {1};
          \node[vertex] (0) at (-2,-0.8) {0};
          \node[emptygraph] (e) at (0,0) {$2 \sim n-1$};
          \draw[-latex] (1) -- (1 -| e.west);
          \draw[latex-] (0) -- (0 -| e.west);
        \end{tikzpicture}\par
        \smallskip\par
        $T_2(n)$
        \bigskip\par
        $A_{2,2}(n): (1,*_{n-2}) \mid (*_{n-2},0)$
      \end{center}
      \caption{$2(n-2)$ secrets}\label{fig:mm:c}
    \end{subfigure}
    \bigskip\par
    \begin{subfigure}[b]{50mm}
      \begin{center}
        \begin{tikzpicture}
          \node[vertex] (0) at (-2,-0.8) {0};
          \node[emptygraph] (e) at (0,0) {$1 \sim n-1$};
          \draw[latex-] (0) -- (0 -| e.west);
        \end{tikzpicture}
        \smallskip\par
        $T_1(n)$
        \bigskip\par
        $A_{2,4}(n): (*_{n-1},0)$
      \end{center}
      \caption{$n-1$ secrets}\label{fig:mm:d}
    \end{subfigure}
    \begin{subfigure}[b]{75mm}
      \begin{center}
        \begin{tikzpicture}
          \node[vertex] (0) at (-2,-0.8) {0};
          \node[emptygraph] (e) at (0,0) {$1 \sim n-1$};
          \draw[-latex] (0) -- (e);
          \draw[latex-] (0) -- (0 -| e.west);
        \end{tikzpicture}
        \smallskip\par
        $T_{12}(n)$
        \bigskip\par
        $A_{2,1}(n): (0,*_{n-1}) \mid (*_{n-1},0)$
      \end{center}
      \caption{$2(n-1)$ secrets}\label{fig:mm:e}
    \end{subfigure}
    \bigskip\par
    \begin{subfigure}[b]{50mm}
      \begin{center}
        \begin{tikzpicture}
          \node[vertex] (0) at (-2,-0.8) {0};
          \node[emptygraph] (e) at (0,0) {$1 \sim n-1$};
          \draw[-latex] (0) -- (0 -| e.west);
          \draw[-latex] (0) edge [loop above] (0);
        \end{tikzpicture}
        \smallskip\par
        $T_{10}(n)$
        \bigskip\par
        $A_{2,5}(n): (0,0) \mid (0,*_{n-1})$
      \end{center}
      \caption{$n$ secrets}\label{fig:mm:f}
    \end{subfigure}
    \begin{subfigure}[b]{75mm}
      \begin{center}
        \begin{tikzpicture}
          \node[vertex] (1) at (-2,0.8) {1};
          \node[vertex] (0) at (-2,-0.8) {0};
          \node[emptygraph] (e) at (0,0) {$2 \sim n-1$};
          \draw[-latex] (1) -- (1 -| e.west);
          \draw[-latex] (1) edge [loop above] (1);
          \draw[latex-] (0) -- (0 -| e.west);
          \draw[-latex] (0) edge [loop above] (0);
        \end{tikzpicture}
        \smallskip\par
        $T_{20}(n)$\
        \bigskip\par
        $A_{2,3}(n): (0,0) \mid (0,*_{n-2}) \mid (1,1) \mid (*_{n-2},1)$
      \end{center}
      \caption{$2(n-1)$ secrets}\label{fig:mm:g}
    \end{subfigure}
  \end{center}
  \caption{Two pegs Mastermind, maset patterns}\label{fig:mm}
\end{figure}

There are 6 maset patterns found by the program.
They corresponds to 6 non-empty (with edges) graphs considered in \cite{ChL04mm}
as shown in Figure \ref{fig:mm}, see also Figure 2 of \cite{ChL04mm}.
Figure \ref{fig:mm} shows the correspondence between our notation $A_{p,i}$ for
the external path length and notation $T_x$, where $T_x$ denotes one of the
function $T$, $T_{1}$, $T_{2}$, $T_{10}$, $T_{12}$, $T_{20}$ considered in
\cite{ChL04mm}.
We use different notation, because we consider the game with the additional
color and because we index functions $A_{p,i}$ in the order they are found by
the program.
From two graphs drawn in Figure 2d of \cite{ChL04mm}, only one is drawn here in
Figure \ref{fig:mm:d}.
The other one corresponds to the maset pattern $(0,*_{n-1})$ which is isomorphic
to the maset pattern $(*_{n-1},0)$.
Similarly in Figure \ref{fig:mm:f},
the other graph drawn in \cite{ChL04mm} corresponds to the maset pattern
$(0,0) \mid (*_{n-1},0)$ which is isomorphic to the maset pattern
$(0,0) \mid (0,*_{n-1})$.

\begin{table}[htb]
  \caption{Two pegs Mastermind, the minimal external path length
           for small number of colors}\label{tab:mm2val}
  \begin{center}
    \begin{tabular}{|c|cccccc|}
\hline
$n$ & $A_{2,0}(n)$ & $A_{2,1}(n)$ & $A_{2,2}(n)$ & $A_{2,3}(n)$ & $A_{2,4}(n)$ & $A_{2,5}(n)$ \\
\hline
2 & 8 & 3 &  &  & 1 & 3 \\
3 & 21 & 7 & 3 & 7 & 3 & 6 \\
4 & 45 & 13 & 7 & 13 & 6 & 9 \\
5 & 81 & 21 & 13 & 21 & 9 & 13 \\
\hline
\end{tabular}

  \end{center}
\end{table}

\begin{table}[htb]
  \caption{Two pegs Mastermind, corresponding equations}\label{tab:mm2eq}
  \begin{center}
    \begin{tabular}{|c|cccccc|}
      \hline
      $j$ & $A_{2,0}$ & $A_{2,1}$ & $A_{2,2}$ & $A_{2,3}$ & $A_{2,4}$ & $A_{2,5}$ \\
      \cline{2-7}
          & $T$       & $T_{12}$  & $T_{2}$   & $T_{20}$  & $T_{1}$   & $T_{10}$  \\
      \hline
      1   & 1         & 3         & 5         & 6         & 4         & 4         \\
      2   & \opt{2}   & \opt{1}   & 4.2       & 4         & 2         & 1         \\
      3   & new       & new       & 2         & 1         & new       & new       \\
      4   &           & 4         & new       & new       & 1         & \opt{2}   \\
      5   &           & 2         & 4.1       & 5         & 5         & 5         \\
      6   &           & new       & \opt{1}   & \opt{2}   & \opt{3}   & 3         \\
      7   &           &           & new       & new       & new       & new       \\
      8   &           &           & 6         & 7         & new       & new       \\
      9   &           &           & 3         & 3         & new       & new       \\
      10  &           &           & new       & new       &           &           \\
      \hline
    \end{tabular}
  \end{center}
\end{table}

Table \ref{tab:mm2val} contains computed values of functions $A_{2,i}$ for small
number of colors.
Comparing this table with Table 2 of \cite{ChL04mm} convinces us that using the
additional color does not help to decrease the number of questions for at most 5
colors, i.e.\ we have that $A_{2,0}(n)=T(n)$, $A_{2,1}(n)=T_{12}(n)$,
$A_{2,2}(n)=T_{2}(n)$, $A_{2,3}(n)=T_{20}(n)$, $A_{2,4}(n)=T_{1}(n)$,
$A_{2,5}(n)=T_{10}(n)$ for $n \le 5$.

There are 47 equations found by the program.
Some equations are not printed in the simplest form.
For example, equation (\ref{eq:mm.2.1.2}) could be simplified as
$A_{2,1}(n) = 2A_{2,4}(n-1) + 2n - 1$.
Such simplification can be implemented, but we do not want to complicate the
program.
Table \ref{tab:mm2eq} shows correspondence between our equations and equations
presented in Table 1 of \cite{ChL04mm}.
Column $A_{2,i}$ and row $j$ contains the case number of function $T_x$ for our
equation (MM.2.$i$.$j$) or the word ``new'' which means that this is a new
equation for a question with the additional color.
We see that our 30 equations for questions without the additional color are
identical as in \cite{ChL04mm}.
It is proved there that for each game state there is always an optimal question.
Equations, which give the optimal solution, are marked in Table \ref{tab:mm2eq}
in bold.
Below, we consider new equations one by one and we show that using the
additional color does not improve the global solution.
We consider equations in reverse index order, because of recursive dependencies.

Equation (\ref{eq:mm.2.5.9}) is the same as equation (\ref{eq:mm.2.5.5}).
Hence, the question $(\ac,1)$ does not give any advantage for the codebreaker
over the question $(1,1)$ in the game state $M_{2,5}$.

Equation (\ref{eq:mm.2.5.8}) gives always more questions than equation
(\ref{eq:mm.2.5.1}).
Hence, the question $(\ac,0)$ does not give any advantage for the codebreaker
over the question $(0,0)$ in the game state $M_{2,5}$.

Equation (\ref{eq:mm.2.5.7}) is the same as equation (\ref{eq:mm.2.5.5}).
Hence, the question $(1,\ac)$ does not give any advantage for the codebreaker
over the question $(1,1)$ in the game state $M_{2,5}$.

Equation (\ref{eq:mm.2.5.3}) has no solution, because the question $(0,\ac)$
does not distinguish secrets in the game state $M_{2,5}$.

We conclude that $A_{2,5}(n)=T_{10}(n)$ for $n \ge 2$ and hence the additional
color does not help in the game state $M_{2,5}$.

Equation (\ref{eq:mm.2.4.9}) is the same as equation (\ref{eq:mm.2.4.5}).
Hence, the question $(\ac,1)$ does not give any advantage for the codebreaker
over the question $(1,1)$ in the game state $M_{2,4}$.

Equation (\ref{eq:mm.2.4.8}) has no solution, because the question $(\ac,0)$
does not distinguish secrets in the game state $M_{2,4}$.

Equation (\ref{eq:mm.2.4.7}) is the same as equation (\ref{eq:mm.2.4.5}).
Hence, the question $(1,\ac)$ does not give any advantage for the codebreaker
over the question $(1,1)$ in the game state $M_{2,4}$.

Equation (\ref{eq:mm.2.4.3}) has no solution, because the question $(0,\ac)$
does not distinguish secrets in the game state $M_{2,4}$.

We conclude that $A_{2,4}(n)=T_{1}(n)$ for $n \ge 2$ and hence the additional
color does not help in the game state $M_{2,4}$.

Equation (\ref{eq:mm.2.3.10}) has the solution $A_{2,3}(n)=n^2+n+C$, where $C$
is a constant depending on the beginning value of recurrence.
If we compare the solution with values of $A_{2,3}(n)$ in Table \ref{tab:mm2val}
and the formula for $T_{20}(n)$ in Table 2 of \cite{ChL04mm} then we see that
the question $(2,\ac)$ is never optimal for the codebreaker in the game state
$M_{2,3}$.

Equation (\ref{eq:mm.2.3.7}) gives always more questions than equation
(\ref{eq:mm.2.3.1}).
Hence, the question $(1,\ac)$ does not give any advantage for the codebreaker
over the question $(0,0)$ in the game state $M_{2,3}$.

Equation (\ref{eq:mm.2.3.4}) can be solved using the know solution
$A_{2,5}(n)=T_{10}(n)$.
By this assumption, we have that $A_{2,3}(n)=2T_{10}(n-1)+2n-2$.
This always gives more questions than $T_{20}(n)$, which can be verified by
Table \ref{tab:mm2val} for $n \le 5$ and by Table 2 of \cite{ChL04mm} for
$n > 5$.
Hence, the question $(0,\ac)$ is never optimal for the codebreaker in the game
state $M_{2,3}$.

We conclude that $A_{2,3}(n)=T_{20}(n)$ for $n \ge 3$ and hence the additional
color does not help in the game state $M_{2,3}$.

Equation (\ref{eq:mm.2.2.10}) has the solution $A_{2,2}(n)=n^2-n+C$, where $C$
is a constant depending on the beginning value of recurrence.
If we compare the solution with values of $A_{2,2}(n)$ in Table \ref{tab:mm2val}
and the formula for $T_{2}(n)$ in Table 2 of \cite{ChL04mm} then we see that
the question $(2,\ac)$ is never optimal for the codebreaker in the game state
$M_{2,2}$.

Equation (\ref{eq:mm.2.2.7}) is the same as equation (\ref{eq:mm.2.2.1}).
Hence, the question $(1,\ac)$ does not give any advantage for the codebreaker
over the question $(0,0)$ in the game state $M_{2,2}$.

Equation (\ref{eq:mm.2.2.4}) is the same as equation (\ref{eq:mm.2.2.1}).
Hence, the question $(0,\ac)$ does not give any advantage for the codebreaker
over the question $(0,0)$ in the game state $M_{2,2}$.

We conclude that $A_{2,2}(n)=T_{2}(n)$ for $n \ge 3$ and hence the additional
color does not help in the game state $M_{2,2}$.

Equation (\ref{eq:mm.2.1.6}) has the solution $A_{2,1}(n)=n^2+n+C$, where $C$
is a constant depending on the beginning value of recurrence.
If we compare the solution with values of $A_{2,1}(n)$ in Table \ref{tab:mm2val}
and the formula for $T_{12}(n)$ in Table 2 of \cite{ChL04mm} then we see that
the question $(1,\ac)$ is never optimal for the codebreaker in the game state
$M_{2,1}$.

Equation (\ref{eq:mm.2.1.3}) can be solved using the know solution
$A_{2,4}(n)=T_{1}(n)$.
By this assumption, we have that $A_{2,1}(n)=2T_{1}(n)+2n-2$.
This always gives more questions than $T_{12}(n)$, which can be verified by
Table \ref{tab:mm2val} for $n \le 5$ and by Table 2 of \cite{ChL04mm} for
$n > 5$.
Hence, the question $(0,\ac)$ is never optimal for the codebreaker in the game
state $M_{2,1}$.

We conclude that $A_{2,1}(n)=T_{12}(n)$ for $n \ge 2$ and hence the additional
color does not help in the game state $M_{2,1}$.

Equation (\ref{eq:mm.2.0.3}) can be solved using the know solutions
$A_{2,4}(n)=T_{1}(n)$ and $A_{2,5}(n)=T_{10}(n)$.
By this assumption, we have that
$A_{2,0}(n)-A_{2,0}(n-1)=T_{1}(n)+T_{10}(n)+n^2$.
This always gives greater value than $T(n)-T(n-1)$, which can be verified by
Table \ref{tab:mm2val} for $n \le 5$ and by Table 2 of \cite{ChL04mm} for
$n > 5$.
Hence, the question $(0,\ac)$ is never optimal for the codebreaker in the game
state $M_{2,0}$.

We conclude that $A_{2,0}(n)=T(n)$ for $n \ge 2$ and hence the additional color
does not help in the game state $M_{2,0}$.
This last conclusion shows that using the additional color does not decrease
the minimal number of questions in the expected case of Mastermind and it ends
the proof of equation~(\ref{eq:mme}).

\subsection{Three pegs}

% 13388 maset-patterns
% 9096599 equations
% dependency cycle not found
% run time without printing
% real    31m44.852s
% user    31m42.871s
% sys     0m0.862s

Authors claim in \cite{ChL04mm} that they technique (called graph-partition) is
not easily extensible to more than two pegs.
Perhaps, they would need to consider multigraphs, which are not easy to draw.
In contrary, our algorithm works for any number of pegs.

For three pegs, there are found 13388 maset patterns and 9096599 equations.
Unfortunately, this is too much to be solved manually.
Even listing all equation is quite big job.
In our experiment, we only counted equations and we did not list them
explicitly.
Further research should concentrate on developing an algorithm or using an
existing software to efficiently solve such big system of equation.
The difficulty lies in proving that indeed for each state of the game, it is
always an optimal question, regardless of the value of $n$, or the choice of
question depends on~$n$.

\section{AB game}\label{sec:ab}

In this section, we present results of computations for AB game.
Next, we analyze the results.
One peg AB game is equivalent to one peg Mastermind, therefore we omit the case
$p=1$.
Compared to Mastermind, the algorithm differs in that we avoid repeating
explicit colors in a clause when generating and extending maset patterns.
We also avoid questions with repeated explicit colors.
When counting the number of secrets represented by the maset pattern with more
than one star in a clause, we must take into account that explicit colors cannot
be repeated, e.g.\ the star $*_{n-t}$ represents $n-t$ colors, but the clause
$(*_{n-t},*_{n-t})$ represents only $(n-t)(n-t-1)$ secrets.

\subsection{Two pegs}

\begin{itemize}
  \item[] $M_{2,0} = (*_{n},*_{n})$
    \begin{itemize}
      \item[] $(0,1)$
        \begin{itemize}
          \item[] $\ans{0}{0} \rightarrow (*_{n-2},*_{n-2})$
          \item[] $\ans{0}{1} \rightarrow (1,*_{n-2}) \mid (*_{n-2},0)$
          \item[] $\ans{1}{0} \rightarrow (0,*_{n-2}) \mid (*_{n-2},1) \Rrightarrow (1,*_{n-2}) \mid (*_{n-2},0)$
          \item[] $\ans{0}{2} \rightarrow (1,0)$
          \item[] $\ans{2}{0} \rightarrow (0,1)$
            \begin{equation} A_{2,0}(n) = A_{2,0}(n-2) + A_{2,1}(n) + A_{2,1}(n) + n(n-1) + 1 \tag{AB.2.0.1} \label{eq:ab.2.0.1} \end{equation}
        \end{itemize}
      \item[] $(0,\ac)$
        \begin{itemize}
          \item[] $\ans{0}{0} \rightarrow (*_{n-1},*_{n-1})$
          \item[] $\ans{0}{1} \rightarrow (*_{n-1},0)$
          \item[] $\ans{1}{0} \rightarrow (0,*_{n-1}) \Rrightarrow (*_{n-1},0)$
          \item[] $\ans{0}{2} \rightarrow \emptyset$
          \item[] $\ans{2}{0} \rightarrow \emptyset$
            \begin{equation} A_{2,0}(n) = A_{2,0}(n-1) + A_{2,2}(n) + A_{2,2}(n) + n(n-1) \tag{AB.2.0.2} \label{eq:ab.2.0.2} \end{equation}
        \end{itemize}
    \end{itemize}
  \item[] $M_{2,1} = (1,*_{n-2}) \mid (*_{n-2},0)$
    \begin{itemize}
      \item[] $(0,1)$
        \begin{itemize}
          \item[] $\ans{0}{0} \rightarrow \emptyset$
          \item[] $\ans{0}{1} \rightarrow (1,*_{n-2}) \mid (*_{n-2},0)$
          \item[] $\ans{1}{0} \rightarrow \emptyset$
          \item[] $\ans{0}{2} \rightarrow \emptyset$
          \item[] $\ans{2}{0} \rightarrow \emptyset$
            \begin{equation} A_{2,1}(n) = A_{2,1}(n) + 2(n-2) \tag{AB.2.1.1} \label{eq:ab.2.1.1} \end{equation}
        \end{itemize}
      \item[] $(0,2)$
        \begin{itemize}
          \item[] $\ans{0}{0} \rightarrow (1,*_{n-3}) \Rightarrow (0,*_{n-3}) \Rrightarrow (*_{n-3},0)$
          \item[] $\ans{0}{1} \rightarrow (*_{n-3},0)$
          \item[] $\ans{1}{0} \rightarrow (1,2)$
          \item[] $\ans{0}{2} \rightarrow (2,0)$
          \item[] $\ans{2}{0} \rightarrow \emptyset$
            \begin{equation} A_{2,1}(n) = A_{2,2}(n-2) + A_{2,2}(n-2) + 2(n-2) + 2 \tag{AB.2.1.2} \label{eq:ab.2.1.2} \end{equation}
        \end{itemize}
      \item[] $(0,\ac)$
        \begin{itemize}
          \item[] $\ans{0}{0} \rightarrow (1,*_{n-2}) \Rightarrow (0,*_{n-2}) \Rrightarrow (*_{n-2},0)$
          \item[] $\ans{0}{1} \rightarrow (*_{n-2},0)$
          \item[] $\ans{1}{0} \rightarrow \emptyset$
          \item[] $\ans{0}{2} \rightarrow \emptyset$
          \item[] $\ans{2}{0} \rightarrow \emptyset$
            \begin{equation} A_{2,1}(n) = A_{2,2}(n-1) + A_{2,2}(n-1) + 2(n-2) \tag{AB.2.1.3} \label{eq:ab.2.1.3} \end{equation}
        \end{itemize}
      \item[] $(1,0)$
        \begin{itemize}
          \item[] $\ans{0}{0} \rightarrow \emptyset$
          \item[] $\ans{0}{1} \rightarrow \emptyset$
          \item[] $\ans{1}{0} \rightarrow (1,*_{n-2}) \mid (*_{n-2},0)$
          \item[] $\ans{0}{2} \rightarrow \emptyset$
          \item[] $\ans{2}{0} \rightarrow \emptyset$
            \begin{equation} A_{2,1}(n) = A_{2,1}(n) + 2(n-2) \tag{AB.2.1.4} \label{eq:ab.2.1.4} \end{equation}
        \end{itemize}
      \item[] $(1,2)$
        \begin{itemize}
          \item[] $\ans{0}{0} \rightarrow (*_{n-3},0)$
          \item[] $\ans{0}{1} \rightarrow (2,0)$
          \item[] $\ans{1}{0} \rightarrow (1,*_{n-3}) \Rightarrow (0,*_{n-3}) \Rrightarrow (*_{n-3},0)$
          \item[] $\ans{0}{2} \rightarrow \emptyset$
          \item[] $\ans{2}{0} \rightarrow (1,2)$
            \begin{equation} A_{2,1}(n) = A_{2,2}(n-2) + A_{2,2}(n-2) + 2(n-2) + 1 \tag{AB.2.1.5} \label{eq:ab.2.1.5} \end{equation}
        \end{itemize}
      \item[] $(1,\ac)$
        \begin{itemize}
          \item[] $\ans{0}{0} \rightarrow (*_{n-2},0)$
          \item[] $\ans{0}{1} \rightarrow \emptyset$
          \item[] $\ans{1}{0} \rightarrow (1,*_{n-2}) \Rightarrow (0,*_{n-2}) \Rrightarrow (*_{n-2},0)$
          \item[] $\ans{0}{2} \rightarrow \emptyset$
          \item[] $\ans{2}{0} \rightarrow \emptyset$
            \begin{equation} A_{2,1}(n) = A_{2,2}(n-1) + A_{2,2}(n-1) + 2(n-2) \tag{AB.2.1.6} \label{eq:ab.2.1.6} \end{equation}
        \end{itemize}
      \item[] $(2,3)$
        \begin{itemize}
          \item[] $\ans{0}{0} \rightarrow (1,*_{n-4}) \mid (*_{n-4},0)$
          \item[] $\ans{0}{1} \rightarrow (1,2) \mid (3,0)$
          \item[] $\ans{1}{0} \rightarrow (1,3) \mid (2,0)$
          \item[] $\ans{0}{2} \rightarrow \emptyset$
          \item[] $\ans{2}{0} \rightarrow \emptyset$
            \begin{equation} A_{2,1}(n) = A_{2,1}(n-2) + 2(n-2) + 6 \tag{AB.2.1.7} \label{eq:ab.2.1.7} \end{equation}
        \end{itemize}
      \item[] $(2,\ac)$
        \begin{itemize}
          \item[] $\ans{0}{0} \rightarrow (1,*_{n-3}) \mid (*_{n-3},0)$
          \item[] $\ans{0}{1} \rightarrow (1,2)$
          \item[] $\ans{1}{0} \rightarrow (2,0)$
          \item[] $\ans{0}{2} \rightarrow \emptyset$
          \item[] $\ans{2}{0} \rightarrow \emptyset$
            \begin{equation} A_{2,1}(n) = A_{2,1}(n-1) + 2(n-2) + 2 \tag{AB.2.1.8} \label{eq:ab.2.1.8} \end{equation}
        \end{itemize}
    \end{itemize}
  \item[] $M_{2,2} = (*_{n-1},0)$
    \begin{itemize}
      \item[] $(0,1)$
        \begin{itemize}
          \item[] $\ans{0}{0} \rightarrow \emptyset$
          \item[] $\ans{0}{1} \rightarrow (*_{n-2},0)$
          \item[] $\ans{1}{0} \rightarrow \emptyset$
          \item[] $\ans{0}{2} \rightarrow (1,0)$
          \item[] $\ans{2}{0} \rightarrow \emptyset$
            \begin{equation} A_{2,2}(n) = A_{2,2}(n-1) + (n-1) + 1 \tag{AB.2.2.1} \label{eq:ab.2.2.1} \end{equation}
        \end{itemize}
      \item[] $(0,\ac)$
        \begin{itemize}
          \item[] $\ans{0}{0} \rightarrow \emptyset$
          \item[] $\ans{0}{1} \rightarrow (*_{n-1},0)$
          \item[] $\ans{1}{0} \rightarrow \emptyset$
          \item[] $\ans{0}{2} \rightarrow \emptyset$
          \item[] $\ans{2}{0} \rightarrow \emptyset$
            \begin{equation} A_{2,2}(n) = A_{2,2}(n) + (n-1) \tag{AB.2.2.2} \label{eq:ab.2.2.2} \end{equation}
        \end{itemize}
      \item[] $(1,0)$
        \begin{itemize}
          \item[] $\ans{0}{0} \rightarrow \emptyset$
          \item[] $\ans{0}{1} \rightarrow \emptyset$
          \item[] $\ans{1}{0} \rightarrow (*_{n-2},0)$
          \item[] $\ans{0}{2} \rightarrow \emptyset$
          \item[] $\ans{2}{0} \rightarrow (1,0)$
            \begin{equation} A_{2,2}(n) = A_{2,2}(n-1) + (n-1) \tag{AB.2.2.3} \label{eq:ab.2.2.3} \end{equation}
        \end{itemize}
      \item[] $(1,2)$
        \begin{itemize}
          \item[] $\ans{0}{0} \rightarrow (*_{n-3},0)$
          \item[] $\ans{0}{1} \rightarrow (2,0)$
          \item[] $\ans{1}{0} \rightarrow (1,0)$
          \item[] $\ans{0}{2} \rightarrow \emptyset$
          \item[] $\ans{2}{0} \rightarrow \emptyset$
            \begin{equation} A_{2,2}(n) = A_{2,2}(n-2) + (n-1) + 2 \tag{AB.2.2.4} \label{eq:ab.2.2.4} \end{equation}
        \end{itemize}
      \item[] $(1,\ac)$
        \begin{itemize}
          \item[] $\ans{0}{0} \rightarrow (*_{n-2},0)$
          \item[] $\ans{0}{1} \rightarrow \emptyset$
          \item[] $\ans{1}{0} \rightarrow (1,0)$
          \item[] $\ans{0}{2} \rightarrow \emptyset$
          \item[] $\ans{2}{0} \rightarrow \emptyset$
            \begin{equation} A_{2,2}(n) = A_{2,2}(n-1) + (n-1) + 1 \tag{AB.2.2.5} \label{eq:ab.2.2.5} \end{equation}
        \end{itemize}
      \item[] $(\ac,0)$
        \begin{itemize}
          \item[] $\ans{0}{0} \rightarrow \emptyset$
          \item[] $\ans{0}{1} \rightarrow \emptyset$
          \item[] $\ans{1}{0} \rightarrow (*_{n-1},0)$
          \item[] $\ans{0}{2} \rightarrow \emptyset$
          \item[] $\ans{2}{0} \rightarrow \emptyset$
            \begin{equation} A_{2,2}(n) = A_{2,2}(n) + (n-1) \tag{AB.2.2.6} \label{eq:ab.2.2.6} \end{equation}
        \end{itemize}
      \item[] $(\ac,1)$
        \begin{itemize}
          \item[] $\ans{0}{0} \rightarrow (*_{n-2},0)$
          \item[] $\ans{0}{1} \rightarrow (1,0)$
          \item[] $\ans{1}{0} \rightarrow \emptyset$
          \item[] $\ans{0}{2} \rightarrow \emptyset$
          \item[] $\ans{2}{0} \rightarrow \emptyset$
            \begin{equation} A_{2,2}(n) = A_{2,2}(n-1) + (n-1) + 1 \tag{AB.2.2.7} \label{eq:ab.2.2.7} \end{equation}
        \end{itemize}
    \end{itemize}
\end{itemize}

% 3 maset-patterns
% 17 equations
% dependency cycle not found

\begin{figure}[htb]
  \tikzstyle{vertex}=[circle, draw, minimum size=8mm, inner sep=0cm]
  \tikzstyle{emptygraph}=[rectangle, draw, minimum width=22mm, minimum height=24mm, inner sep=0cm]
  \tikzstyle{fullgraph}=[ellipse, draw, minimum width=22mm, minimum height=28mm, inner sep=0cm]
  \begin{center}
    \begin{subfigure}[b]{50mm}
      \begin{center}
        \begin{tikzpicture}
          \node[emptygraph] at (0,0) {$0 \sim n-1$};
        \end{tikzpicture}
        \smallskip\par
        \null\par
        \bigskip\par
        $\emptyset$
      \end{center}
      \caption{no secrets}\label{fig:ab:a}
    \end{subfigure}
    \begin{subfigure}[b]{50mm}
      \begin{center}
        \begin{tikzpicture}
          \node[fullgraph] at (0,0) {$0 \sim n-1$};
        \end{tikzpicture}
        \smallskip\par
        $T(n)$\par
        \bigskip\par
        $A_{2,0}(n): (*_{n},*_{n})$\par
      \end{center}
      \caption{$n(n-1)$ secrets}\label{fig:ab:b}
    \end{subfigure}
    \bigskip\par
    \begin{subfigure}[b]{50mm}
      \begin{center}
        \begin{tikzpicture}
          \node[vertex] (0) at (-2,-0.8) {0};
          \node[emptygraph] (e) at (0,0) {$1 \sim n-1$};
          \draw[latex-] (0) -- (0 -| e.west);
        \end{tikzpicture}
        \smallskip\par
        $T_1(n)$
        \bigskip\par
        $A_{2,2}(n): (*_{n-1},0)$
      \end{center}
      \caption{$n-1$ secrets}\label{fig:ab:c}
    \end{subfigure}
    \begin{subfigure}[b]{50mm}
      \begin{center}
        \begin{tikzpicture}
          \node[vertex] (1) at (-2,0.8) {1};
          \node[vertex] (0) at (-2,-0.8) {0};
          \node[emptygraph] (e) at (0,0) {$2 \sim n-1$};
          \draw[-latex] (1) -- (1 -| e.west);
          \draw[latex-] (0) -- (0 -| e.west);
        \end{tikzpicture}\par
        \smallskip\par
        $T_2(n)$
        \bigskip\par
        $A_{2,1}(n): (1,*_{n-2}) \mid (*_{n-2},0)$
      \end{center}
      \caption{$2(n-2)$ secrets}\label{fig:ab:d}
    \end{subfigure}
  \end{center}
  \caption{Two pegs AB game, maset patterns}\label{fig:ab}
\end{figure}
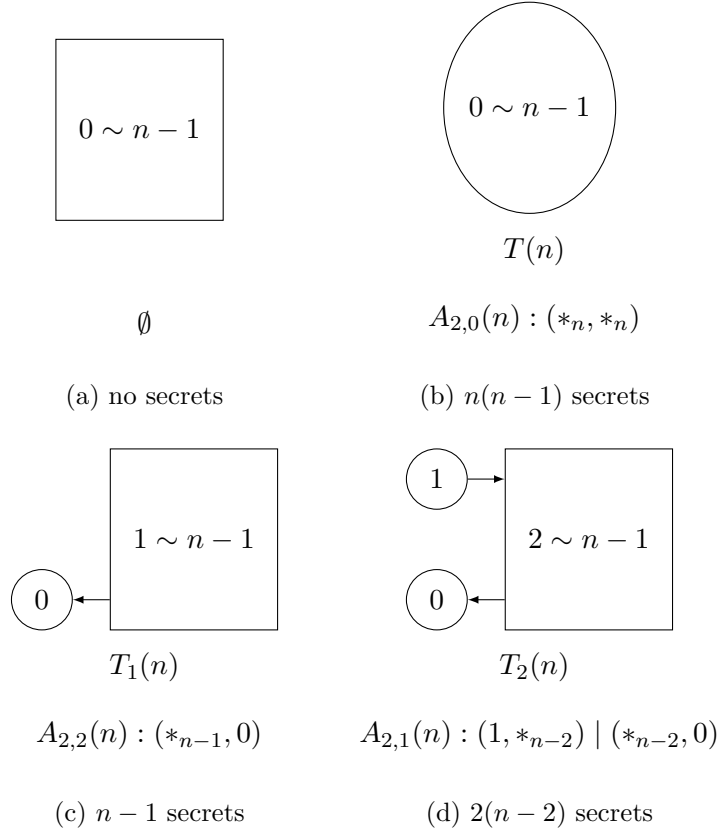

There are 3 maset patterns found by the program.
They corresponds to 3 non-empty  (with edges) graphs considered in
\cite{ChL04ab} as shown in Figure \ref{fig:ab}, see also Figure 5 of
\cite{ChL04ab}.
Figure \ref{fig:ab} shows the correspondence between our notation $A_{p,i}$ for
the external path length and notation $T_x$, where $T_x$ denotes one of the
function $T$, $T_1$, $T_2$ considered in \cite{ChL04ab}.
As previously, we use different notation, because we consider the game with the
additional color and because we index functions $A_{p,i}$ in the order they are
found by the program.
From two graphs drawn in Figure 5c of \cite{ChL04ab}, only one is drawn here in
Figure \ref{fig:ab:c}.
The other one corresponds to the maset pattern $(0,*_{n-1})$ which is isomorphic
to the maset pattern $(*_{n-1},0)$.

\begin{table}[htb]
  \caption{Two pegs AB game, the minimal external path length for
           small number of colors}\label{tab:ab2val}
  \begin{center}
    \begin{tabular}{|c|ccc|}
\hline
$n$ & $A_{2,0}(n)$ & $A_{2,1}(n)$ & $A_{2,2}(n)$ \\
\hline
2 & 3 &  & 1 \\
3 & 13 & 3 & 3 \\
4 & 30 & 7 & 6 \\
5 & 60 & 13 & 9 \\
\hline
\end{tabular}

  \end{center}
\end{table}

\begin{table}[htb]
  \caption{Two pegs AB game, corresponding cases}\label{tab:ab2eq}
  \begin{center}
    \begin{tabular}{|c|ccc|}
      \hline
      $j$ & $A_{2,0}$ & $A_{2,1}$ & $A_{2,2}$   \\
      \cline{2-4}
          & $T$       & $T_2$     & $T_1$       \\
      \hline
      1   & \opt{0}   & 4         & (ii)        \\
      2   & new       & 2         & new         \\
      3   &           & new       & (i)         \\
      4   &           & 4         & \opt{(iii)} \\
      5   &           & \opt{1}   & new         \\
      6   &           & new       & new         \\
      7   &           & 3         & new         \\
      8   &           & new       &             \\
      \hline
    \end{tabular}
  \end{center}
\end{table}

Table \ref{tab:ab2val} contains computed values of functions $A_{2,i}$ for small
number of colors.
Comparing this table with Table 1 of \cite{ChL04ab} convinces us that using the
additional color does not help to decrease the number of questions for at most 5
colors, i.e.\ we have that $A_{2,0}(n)=T(n)$, $A_{2,1}(n)=T_2(n)$,
$A_{2,2}(n)=T_1(n)$ for $n \le 5$.

There are 17 equations found by the program.
Some equations are not printed in the simplest form.
As previously, we do not simplify equations, because we do not want to
complicate the program.
Table \ref{tab:ab2eq} shows correspondence between our equations and equations
considered in \cite{ChL04ab}.
Column $A_{2,i}$ and row $j$ contains the case number of function $T_x$ for our
equation (AB.2.$i$.$j$) or the word ``new'' which means that this is a new
equation for a question with the additional color.
Zero in the first column means that there is only one case in \cite{ChL04ab},
because there is only one question considered.
Case 4 is listed two times in the second column, because our two questions are
considered in \cite{ChL04ab} as one case leading to the same equation.
We see that our 9 equations for questions without the additional color are
identical as in \cite{ChL04ab}.
It is proved there that for each game state there is always an optimal question.
Cases, which give the optimal solution, are marked in Table \ref{tab:ab2eq} in
bold.
Below, we consider new equations one by one and we show that using the
additional color does not improve the global solution.

Equation (\ref{eq:ab.2.2.7}) is the same as equation (\ref{eq:ab.2.2.1}).
Hence, the question $(\ac,1)$ does not give any advantage for the codebreaker
over the question $(0,1)$ in the game state $M_{2,2}$.

Equation (\ref{eq:ab.2.2.6}) has no solution, because the question $(\ac,0)$
does not distinguish secrets in the game state $M_{2,2}$.

Equation (\ref{eq:ab.2.2.5}) is the same as equation (\ref{eq:ab.2.2.1}).
Hence, the question $(1,\ac)$ does not give any advantage for the codebreaker
over the question $(0,1)$ in the game state $M_{2,2}$.

Equation (\ref{eq:ab.2.2.2}) has no solution, because the question $(0,\ac)$
does not distinguish secrets in the game state $M_{2,2}$.

We conclude that $A_{2,2}(n)=T_1(n)$ for $n \ge 2$ and hence the additional
color does not help in the game state $M_{2,2}$.

Equation (\ref{eq:ab.2.1.8}) has the solution $A_{2,1}(n)=n^2-n+C$, where $C$
is a constant depending on the beginning value of recurrence.
If we compare the solution with values of $A_{2,1}(n)$ in Table \ref{tab:ab2val}
and the formula for $T_2(n)$ in Table 1 of \cite{ChL04mm} then we see that the
question $(2,\ac)$ is never optimal for the codebreaker in the game state
$M_{2,1}$.

Equation (\ref{eq:ab.2.1.6}) is the same as equation (\ref{eq:mm.2.1.3}).
Hence, the question $(1,\ac)$ does not give any advantage for the codebreaker
over the question $(0,\ac)$ in the game state $M_{2,1}$.

Equation (\ref{eq:ab.2.1.3}) can be solved using the know solution
$A_{2,2}(n)=T_1(n)$.
By this assumption, we have that $A_{2,1}(n)=2T_1(n-1)+2n-4$.
This always gives more questions than $T_2(n)$, which can be verified by Table
\ref{tab:ab2val} for $n \le 5$ and by Table 1 of \cite{ChL04ab} for $n > 5$.
Hence, the question $(0,\ac)$ is never optimal for the codebreaker in the game
state $M_{2,1}$.

We conclude that $A_{2,1}(n)=T_2(n)$ for $n \ge 3$ and hence the additional
color does not help in the game state $M_{2,1}$.

Equation (\ref{eq:ab.2.0.2}) can be solved using the know solution
$A_{2,2}(n)=T_1(n)$.
By this assumption, we have that $A_{2,0}(n)-A_{2,0}(n-1)=2T_1(n)+n^2-n$.
This always gives greater value than $T(n)-T(n-1)$, which can be verified by
Table \ref{tab:ab2val} for $n \le 5$ and by Table 1 of \cite{ChL04ab} for
$n > 5$.
Hence, the question $(0,\ac)$ is never optimal for the codebreaker in the game
state~$M_{2,0}$.

We conclude that $A_{2,0}(n)=T(n)$ for $n \ge 2$ and hence the additional color
does not help in the game state $M_{2,0}$.
This last conclusion shows that using the additional color does not decrease
the minimal number of questions in the expected case of AB game and it ends
the proof of equation~(\ref{eq:abe}).

\subsection{Three pegs}

% 7496 maset-patterns
% 4188421 equations
% dependency cycle not found
% run time without printing
% real    13m58.537s
% user    13m57.919s
% sys     0m0.169s

For three pegs, there are found 7496 maset patterns and 4188421 equations.
As previously, this is too much to be solved in a reasonable time.

\end{document}